\documentclass[11pt,a4paper]{amsart}
\usepackage[usenames,dvipsnames]{xcolor}
\usepackage{amssymb,amsmath,epsfig,graphics,mathrsfs,enumerate,verbatim}
\usepackage[pagebackref,colorlinks=true,linkcolor=blue,citecolor=blue]{hyperref}

\usepackage{fancyhdr}
\hypersetup{
 colorlinks=true,
 urlcolor=blue,
 linkcolor=blue,
 citecolor=red,
 bookmarksopen=true
}

\usepackage{amsmath}
\usepackage{amsfonts}
\usepackage{amssymb}
\usepackage{amsthm}
\usepackage{epsfig,graphics,mathrsfs}
\usepackage{graphicx}
\usepackage{tikz}
\def \phi {\varphi}

\def \R {\mathbb{R}}

\newcommand{\Om}{\Omega}

\newcommand{\eps}{\ve}

\newcommand{\bG}{\mathbb {G}}

\newcommand{\la}{\lambda}

\numberwithin{equation}{section}

\newcommand{\beq}{\begin{equation}}
\newcommand{\bea}[1]{\begin{array}{#1} }
\newcommand{\eeq}{ \end{equation}}
\newcommand{\ea}{ \end{array}}

\newcommand{\ve}{\varepsilon}

\newcommand{\x}{\nabla_\mathscr X}

\newtheorem{theorem}{Theorem}[section]

\newtheorem{proposition}[theorem]{Proposition}

\newtheorem{remark}[theorem]{Remark}

\newtheorem{definition}[theorem]{Definition}

\numberwithin{equation}{section}

\usepackage{amscd,amstext}
\usepackage{microtype}
\usepackage{comment}
\usepackage{mathtools}

\providecommand{\rr}{{\mathbb{R}}}

\providecommand{\nn}{{\mathbb{N}}}

\providecommand{\escpr}[1]{\langle#1\rangle}

\providecommand{\frak}{\mathfrak}

\title[]{Exterior $p$-capacitary rigidity on groups of Heisenberg type}


\author{Federica Dragoni}
\address{School of Mathematics\\
Cardiff University\\
Cardiff CF24 4AG
WALES} \email[Federica Dragoni]{DragoniF@cardiff.ac.uk}

\author{Nicola Garofalo}
\address{School of Mathematical and Statistical Sciences\\ Arizona State University}\email[Nicola Garofalo]{nicola.garofalo@asu.edu}

\author{Gianmarco Giovannardi}
\address{Dipartimento di Matematica e Informatica "U. Dini" \\ Universit\`a degli Studi di Firenze\\ 50134 Florence, Italy}
\email[Gianmarco Giovannardi]{gianmarco.giovannardi@unifi.it}

\author{Paolo Salani}
\address{Dipartimento di Matematica e Informatica "U. Dini" \\ Universit\`a degli Studi di Firenze\\
50134 Florence,
ITALY}\email[Paolo Salani]{paolo.salani@unifi.it}

\author{Dimiter Vassilev}
\address{ Department of Mathematics and Statistics\\
University of New Mexico\\
Albuquerque, New Mexico, 87131-0001\\
}
\email[Dimiter Vassilev]{vassilev@unm.edu}

\date{}

\begin{document}

\begin{abstract}
We study the exterior Bernoulli overdetermined problem for sub-elliptic operators.
Let $\bG$ be a group of Heisenberg type with Property (H), and let $\Omega \subset \bG$
be a domain enjoying a partial symmetry property and containing a
Kor\'{a}nyi-Kaplan ball $\overline{B_r}$ centered at the origin.
We consider a function $u$ that is $p$-harmonic with respect to the $p$-sub-Laplacian $\Delta_{\bG,p}$ in the annular region
$\Omega \setminus \overline{B}_r$, satisfying the boundary conditions
$
    u = 0 \;\text{ on }\; \partial\Omega,
    \quad
    u = 1 \;\text{ on }\; \partial B_r.$
We show that if $p$ is even and $u$ satisfies a mild overdetermined condition on the
sub-Riemannian norm of its gradient on the outer boundary $\partial\Omega$,
then $\Omega$ must itself be a Kor\'{a}nyi--Kaplan ball centered at the origin, and $u$ is the corresponding normalized capacitary potential; equivalently, away from the pole it is an affine normalization of the fundamental solution of $\Delta_{\bG,p}$.
\end{abstract}

\maketitle
\tableofcontents

\section{Introduction}

Let $1<p<\infty$, let $0<r<R$, and consider the capacitary potential of the
Euclidean ring $B(0,R)\setminus\overline{B(0,r)}$.  It is the unique solution of
\[
\begin{cases}
\Delta_p u=0&\text{in } B(0,R)\setminus\overline{B(0,r)},\\
u=1&\text{on }\partial B(0,r),\\
u=0&\text{on }\partial B(0,R).
\end{cases}
\]
The solution is radial, and the modulus of its gradient is constant on each
boundary component.  The corresponding overdetermined problem asks whether a
domain with the same property must be a ball.  This question belongs to the
circle of ideas originating with Serrin's theorem \cite{Serrin1971}.  For
exterior Bernoulli problems on ring-shaped domains we refer to
\cite{Al92,Greco2001,HP,HPP}. For a recent extension of the technique of touching balls to the Riemannian setting, see~\cite{GLS26}.

A related, but geometrically different, problem is the capacitary problem in
the unbounded exterior of a bounded set. Reichel \cite{Reichel96} proved
radial symmetry for a broad class of nonlinear overdetermined equations on
smooth exterior domains by the method of moving hyperplanes. For the
$p$-capacitary problem, Garofalo and Sartori \cite{GS99} gave a proof based on
a maximum principle, blow-up and a priori estimates. Their result does not
assume smoothness of the compact conductor, but it assumes starlikeness. This assumption was later removed in \cite{Po19}, in $\rr^d$, $1<p\le d$. In groups of Heisenberg type with Property $(H)$, the $p$-capacitary problem in the unbounded exterior of a partially symmetric set is treated in \cite[Theorem 1.10]{GV23}, by the reduction used here. These
papers concern an unbounded exterior domain, whereas the problem studied here
is posed in a bounded annular domain with a prescribed inner ball and an
unknown outer boundary.

In this paper we prove rigidity results for the corresponding exterior
problem for the horizontal $p$-Laplacian on a group of Heisenberg type.  One key aspect of our result is that the
outer boundary is not assumed to be smooth.  Our first theorem uses an
asymptotic condition on the horizontal gradient in an interior collar of the
boundary.  A second theorem requires only two inequalities, at points where
the boundary has minimum and maximum distance from the origin.  In both
cases the conclusion is that the domain is a Kor\'anyi--Kaplan ball and the
solution is its normalized $p$-capacitary potential.

Let $\bG$ be a group of Heisenberg type, identified with
$\rr^m\times\rr^k$, and whenever convenient we write a generic point $x\in \bG$ in its logarithmic coordinates $(z,\sigma)$.  We use
the Kor\'anyi-Kaplan gauge
\begin{equation}\label{g}
 |x|_{\bG} = |(z,\sigma)|_{\bG}
 =\bigl(|z|^4+16|\sigma|^2\bigr)^{1/4},
\end{equation}
see \cite{Cy81,Ka},
and denote the gauge ball by
\[
B_R (0)=\{(z,\sigma) \in \bG \mid \  \left( | z |^4+ 16 |\sigma|^2\right) < R^4 \}.
\]

If
$Q=m+2k$ is the homogeneous dimension, the radial solutions of the horizontal
$p$-Laplace equation are generated by
\[
 \Gamma_p(x)=
 \begin{cases}
 |x|_{\bG}^{\frac{p-Q}{p-1}},&p\ne Q,\\
 \log |x|_{\bG},&p=Q.
 \end{cases}
\]
Consequently, the capacitary potential of
$B_R(0)\setminus\overline{B_r(0)}$ is
\begin{equation}\label{eq:u}
u(x)=
\begin{cases}
\displaystyle
\frac{|x|_{\bG}^{\frac{p-Q}{p-1}}
      -R^{\frac{p-Q}{p-1}}}
     {r^{\frac{p-Q}{p-1}}-R^{\frac{p-Q}{p-1}}},&p\ne Q,\\[7pt]
\displaystyle
\frac{\log\bigl(|x|_{\bG}/R\bigr)}
     {\log(r/R)},&p=Q.
\end{cases}
\end{equation}
For the nonlinear capacitary potential and capacity of Kor\'anyi--Kaplan rings
see \cite{CDG96}.  On the outer boundary the function in \eqref{eq:u}
satisfies
\[
 |\nabla_{\bG}u(z,\sigma)|=c|z|
\]
for a constant $c>0$.  Thus the natural overdetermined datum is not constant, and
it vanishes on the characteristic set, where $z=0$.

We consider a connected bounded open set $\Omega\subset\bG$ containing
$\overline{B_r(0)}$ and a weak solution of
\begin{equation}\label{eq:DP}
\begin{cases}
\Delta_{\bG,p}u=0
   &\text{in }\Omega\setminus\overline{B_r(0)},\\
u=0&\text{on }\partial\Omega,\\
u=1&\text{on }\overline{B_r(0)}.
\end{cases}
\end{equation}
The domain is assumed to have the partial symmetry described in
Definition~\ref{df:partialsy}, and $\bG$ is assumed to satisfy Property
$(H)$ of Definition~\ref{D:trans}.  By uniqueness and Property $(H)$, the
solution has the form
\[
 u(z,\sigma)=g\left(\frac{|z|^2}{4},\sigma\right).
\]

In the setting of groups of Heisenberg type we assume  that $p$ is an even integer.  Set
\[
 n=\frac{m+p}{2}+k,
 \qquad
 \Omega^{\star\star}
 =\left\{(w,\sigma)\in
 \rr^{\frac{m+p}{2}}\times\rr^k:
 (|w|,\sigma)\in\Omega^\star_+\right\},
\]
where $\Omega^\star_+$ is defined in Definition \ref{df:partialsy}
 and define
\[
 G(w,\sigma)=g(|w|,\sigma).
\]
The choice of the dimension is dictated by the change of variables
$|z|^2/4=|w|$.  This dimension-raising reduction is the one introduced in \cite[Section 6]{GV23}, where it is used for a Serrin-type problem and for the $p$-capacitary problem in the unbounded exterior of a set, see \cite[Theorems 1.7 and 1.10]{GV23}. Here it is carried out at the level of weak solutions, and applied to a bounded ring with a free outer boundary.
Proposition~\ref{pr:G-weak} shows that $G$ is a weak solution of the standard
Euclidean $p$-Laplace equation in
\[
 \Omega^{\star\star}\setminus
 \overline{B^{\star\star}(0,\bar r)},
 \qquad \bar r=\frac{r^2}{4}.
\]
The relation between the two gradients is
\begin{equation}\label{eq:intro-gradient-relation}
 |\nabla_{\bG}u(z,\sigma)|
 =\frac{|z|}{2}|\nabla G(w,\sigma)|.
\end{equation}
More precisely, set $N=\frac{m+p}{2}\in\nn$, and 
for every $w\in\mathbb R^N$ with $|w|=|z|^2/4$,
\[
|\nabla_{\bG}u(z,\sigma)|=\frac{|z|}{2}|\nabla G(w,\sigma)|.
\]
Since $G$ is radial in $w$, the right-hand side is independent of the choice of $w$.
No regularity of $u$ beyond $W^{1,p}_{\bG}$ is needed in the reduction. Once it has been made, the classical Euclidean regularity theory \cite{DiB83,Lewis,To84} applies to $G$, and the continuity of $X_ju$ follows, since $u$ is the composition of $G$ with a smooth map. Regularity of weak solutions with special symmetries can also be obtained directly on the group, see \cite{G09grad} and Section~\ref{sc:BPmeasure}.

Our first result assumes that, for some $c>0$, the quotient
$|\nabla_{\bG}u|/|z|$ tends uniformly to $c$ in a collar of the outer
boundary; the precise hypothesis is \eqref{eq:sRocnb}.  In view of
\eqref{eq:intro-gradient-relation}, the transformed function satisfies
$|\nabla G|\to2c$ near $\partial\Omega^{\star\star}$.  Theorem
\ref{th:euclideointerno} proves that the Euclidean domain
$\Omega^{\star\star}$ is a ball.  Its proof uses regular level sets of $G$ and
comparison with the radial $p$-capacitary potentials of the largest inscribed
and the smallest circumscribed concentric balls.  Theorem
\ref{th:subapprosimatointerno} then gives
\[
 \Omega=B_R(0)
\]
and identifies $u$ with the function in \eqref{eq:u}.  No smoothness of
$\partial\Omega$ is assumed.

We also prove a version in which no condition is imposed along the whole
outer boundary.  Let $x_1$ and $x_2$ be points where
$\partial\Omega^{\star\star}$ has respectively minimum and maximum Euclidean
distance from the origin.  Theorem~\ref{th:eucsharp} assumes only the two
opposite inequalities
\[
 \limsup_{\substack{x\to x_1\\x\in\Omega^{\star\star}}}
 \frac{G(x)}{|x-x_1|}\le c,
 \qquad
 \limsup_{\substack{x\to x_2\\x\in\Omega^{\star\star}}}
 \frac{G(x)}{|x-x_2|}\ge c.
\]
Comparison with the two radial potentials forces the inner and outer radii to
coincide.  Theorem~\ref{th:subapplicationtosharp} gives the corresponding
statement on $\bG$.  The hypotheses there are written directly in terms
of $u$ and its horizontal gradient at the two contact points.

Compared with \cite{GV23}, the new points are the following. The reduction to the Euclidean problem is carried out for weak solutions in $W^{1,p}_{\bG,0}$, with no a priori regularity of $u$ (Proposition~\ref{pr:G-weak}). No regularity of $\partial\Omega$ is assumed. The two-point version, Theorem~\ref{th:subapplicationtosharp}, is sharp for the model: for $\Omega=B_R(0)$, \eqref{eq:limsuples} holds with equality at every point of $\partial\Omega$, \eqref{subcond:exterior_sharp} with equality at the points with $\sigma_2=0$, and \eqref{subcond:exterior2} with equality at the poles.

The paper is organized as follows.  Section~2 contains the notation and the
facts about groups of Heisenberg type that are used later.  Section~3 records
the (formal) variational formulation.
Section~4 contains the reduction to the Euclidean problem and the first sub-Riemannian rigidity theorem.  Section~5 contains
the Euclidean comparison result used in that proof.  Section~6 proves the
two-point Euclidean theorem and its application to groups of Heisenberg type.
\section{Preliminaries}\label{sc:preliminaries}

We use the notation of \cite{GV23} and recall only the facts needed below.

\begin{definition}
A Heisenberg-type algebra is a finite-dimensional real Lie algebra $\mathfrak{g}$ which can be endowed with an inner product $\escpr{\cdot,\cdot}$ such that
\[
[\frak{z}^{\perp},\frak{z}^{\perp}]= \frak{z},
\]
where $\frak{z}$ is the center of $\mathfrak{g}$ and moreover, for every fixed $\sigma \in \frak{z}$  the map
\begin{equation}\label{e:Jz def}
J_{\sigma}: \frak{z}^{\perp} \to \frak{z}^{\perp}
\end{equation}
defined by
\begin{equation}
\label{eq:J}
\escpr{J_{\sigma}(v),w}=\escpr{\sigma ,[v,w]} \quad \forall v,w \in \frak{z}^{\perp}
\end{equation}
is orthogonal whenever $|\sigma|=1$. Equivalently,
\[
J_\sigma^2=-|\sigma|^2I.
\]
Thus, for every $\sigma\ne0$, the map $J_\sigma/|\sigma|$ is an orthogonal
complex structure on $\mathfrak z^\perp$.
\end{definition}
A group of Heisenberg type $\bG$ ($H$-type group, in short) is a connected and simply
connected Lie group whose Lie algebra is an $H$-type algebra.

Put $\mathfrak h=\mathfrak z^\perp$, $m=\dim\mathfrak h$ and
$k=\dim\mathfrak z$.  We identify $\bG$ with its Lie algebra by means of the
exponential map and write a point as $x=(z,\sigma)\in\mathfrak h\oplus
\mathfrak z$.  In these logarithmic coordinates the group law is
\[
(z,\sigma)\circ(z',\sigma')
=\left(z+z',\sigma+\sigma'+\frac12[z,z']\right).
\]
Since $J_\sigma$ is a complex
structure when $|\sigma|=1$, the integer $m$ is even.
Fix orthonormal bases $\{e_1,\ldots,e_m\}$ of $\mathfrak h$ and
$\{\varepsilon_1,\ldots,\varepsilon_k\}$ of $\mathfrak z$.  If $L_x$ denotes
left translation, set
\[
X_j=dL_x(e_j),\qquad T_\ell=dL_x(\varepsilon_\ell).
\]
Thus $X_1,\ldots,X_m$ form an orthonormal frame of the horizontal bundle and
$T_\ell=\partial_{\sigma_\ell}$.  In logarithmic coordinates,
\begin{equation}\label{eq:Xj}
X_j=\frac{\partial}{\partial z_j}
+\frac12\sum_{\ell=1}^k
\langle J_{\varepsilon_\ell}z,e_j\rangle
\frac{\partial}{\partial\sigma_\ell}.
\end{equation}
The defining relation \eqref{eq:J} and the identity
$J_\sigma^2=-|\sigma|^2I$ give
\begin{equation}\label{eq:Jorth}
\langle J_{\varepsilon_\ell}z,J_{\varepsilon_s}z\rangle
=|z|^2\delta_{\ell s},
\qquad
\langle J_{\varepsilon_\ell}z,z\rangle=0.
\end{equation}

For a smooth function $u$ we write
\[
\nabla_{\bG}u=(X_1u,\ldots,X_mu),
\qquad
\Delta_{\bG}u=\sum_{j=1}^mX_j^2u.
\]
For $1<p<\infty$ the horizontal $p$-Laplacian is
\[
\Delta_{\bG,p}u
=\sum_{j=1}^mX_j\bigl(|\nabla_{\bG}u|^{p-2}X_ju\bigr).
\]

\section{The variational setting}
This section provides a heuristic motivation for the subsequent ones. The computations involving the shape derivative are purely formal and the reader should be aware that nothing in Sections 4--6 depends on it.

Through the exponential map we identify $\bG$ with its Lie algebra and write
$dx$ for Lebesgue measure in these coordinates. This is a bi-invariant Haar
measure on $\bG$, see \cite{CG}.
Since such measure interacts with the anisotropic dilations
\[
\delta_\lambda(z,\sigma)=(\lambda z,\lambda^2\sigma)
\]
as follows: $d(\delta_\la(x)) = \la^{m+2k} dx$, the homogeneous dimension is $Q=m+2k$.

Let $1<p<+\infty$.
We consider a bounded open set $\Omega \subset \bG$
with $\overline{B_r(0)} \subset \Omega$ and \emph{Lipschitz boundary} $\partial\Omega$,
so that the outward unit Euclidean normal $\nu$ exists
$d\sigma$-almost everywhere on $\partial\Omega$ by Rademacher's theorem, where
$d\sigma$ denotes Euclidean $(m+k-1)$-dimensional surface measure. For an open set $A \subset \bG$ we set
\[
W_{\bG}^{1,p}(A)= \left\{u \in L^p(A,dx) \mid \ \int_{A} |\nabla_{\bG} u |^p\,dx < +\infty  \right\}
\]
and
\[
\|u\|_{W_{\bG}^{1,p}(A)}= \left(\int_{A} | u |^p\,dx + \int_{A} |\nabla_{\bG} u |^p\,dx\right)^{1/p},
\]
then we set $W_{\bG, 0}^{1,p}(A) =\overline{C_0^{\infty}(A)}^{W_{\bG}^{1,p}}$.
We define the energy functional for the pair $(u, \Omega)$, where
$u \in W_{\bG}^{1,p}(\Omega \setminus \overline{B_r(0)})$ satisfies $u = 1$ on $\partial B_r(0)$
and $u \geq 0$, by
\begin{equation}
\label{eq:functional}
F_{\bG}(u, \Omega)
= \frac{1}{p}\int_{\Omega \setminus \overline{B_r(0)}} |\nabla_{\bG} u|^p \, dx
+ \frac{(p-1)c^p}{p} \int_{\Omega \setminus \overline{B_r(0)}} |z|^p \, dx,
\end{equation}
where $c > 0$ is a given constant and $|z|$ denotes the Euclidean norm
of the horizontal component of $(z, \sigma) \in \bG$.
The two terms compete: the first penalizes the horizontal elastic energy of $u$,
while the second penalizes the weighted volume of the annular region
$\Omega \setminus B_r (0)$ with weight $|z|^p$ reflecting the sub-Riemannian geometry of $\bG$.
Fix $\Omega$ and let $u_t = u + t\varphi$ with $\varphi \in C^\infty_0(\Omega \setminus \overline{B_r(0)})$.
Since the second term in \eqref{eq:functional} does not depend on $u$, only the Dirichlet term contributes
\begin{equation*}
\label{eq:weaksol}
\frac{d}{dt}\bigg|_{t=0} F_{\bG}(u_t, \Omega)
= \int_{\Omega \setminus \overline{B_r(0)}} |\nabla_{\bG} u|^{p-2} \langle \nabla_{\bG} u,\, \nabla_{\bG} \varphi \rangle \, dx = 0.
\end{equation*}
\begin{definition}
We say that $u \in W_{\bG,0}^{1,p}(\Omega)$, with $u=1$ a.e. in $B_r(0)$, is a weak solution of the homogeneous $p$-sub-Laplacian equation in $\Omega \setminus \overline{B_r}(0)$ with the boundary data $u=0$ on $\partial\Omega$, $u=1$ on $\partial B_r(0)$, if
\[
\int_{\Omega \setminus \overline{B_r(0)}} |\nabla_{\bG} u|^{p-2} \langle \nabla_{\bG} u,\, \nabla_{\bG} \varphi \rangle \, dx = 0
\]
for all $\varphi \in C^\infty_0(\Omega \setminus \overline{B_r(0)} )$.
\end{definition}
Integrating by parts via the sub-Riemannian Green formula, and using the fact
that the horizontal vector fields $X_j$ are skew-adjoint with respect to $dx$, yields
\begin{equation*}
- \int_{\Omega \setminus \overline{B_r(0)}} \Delta_{\bG,p} u \cdot \varphi \, dx = 0
\qquad \forall\, \varphi \in C^\infty_0(\Omega \setminus \overline{B_r(0)}).
\end{equation*}
By the fundamental lemma of the calculus of variations we obtain the first
Euler--Lagrange equation:
\begin{equation}
\label{eq:EL1}
\Delta_{\bG,p} u = 0 \qquad \text{in } \Omega \setminus \overline{B_r(0)}.
\end{equation}

Let $\Phi_t = \mathrm{Id} + tV$ be a family of diffeomorphisms of $\bG$
generated by a vector field $V \in C^\infty_c(\bG \setminus \overline{B_r(0)})$,
with $\Omega_t = \Phi_t(\Omega)$. Formally, i.e. assuming that $u$ is $C^1$ up to $\partial\Omega$,
an easy variational computation shows
\begin{equation*}
\frac{d}{dt}\bigg|_{t=0}
\frac1p\int_{\Omega_t \setminus \overline{B_r(0)} }|\nabla_{\bG} u|^p \, dx
= -\frac{p-1}{p}\int_{\partial\Omega} |\nabla_{\bG} u|^p \langle V, \nu \rangle \, d\sigma,
\end{equation*}

\begin{equation*}
\frac{d}{dt}\bigg|_{t=0}
\frac{(p-1)c^p}{p} \int_{\Omega_t \setminus \overline{B_r(0)}} |z|^p \, dx
= \frac{(p-1)c^p}{p} \int_{\partial\Omega} |z|^p \langle V, \nu \rangle \, d\sigma,
\end{equation*}
where $\nu$ denotes the outward Euclidean unit normal to $\partial\Omega$,
which exists $d\sigma$-a.e.\ since $\partial\Omega$ is Lipschitz.
Stationarity of $F_{\bG}$ with respect to all such domain perturbations requires
\begin{equation*}
\int_{\partial\Omega}
\frac{p-1}{p}\Bigl( c^p |z|^p- |\nabla_{\bG} u|^p  \Bigr)
\langle V, \nu \rangle \, d\sigma = 0
\qquad \forall\, V \in C^\infty_c(\bG \setminus \overline{B_r(0)}),
\end{equation*}
which by the fundamental lemma yields the second Euler--Lagrange condition:
\begin{equation}
\label{eq:EL2}
|\nabla_{\bG} u(z,\sigma)|^p = c^p |z|^p
\qquad \text{on } \partial\Omega,
\end{equation}
or equivalently $|\nabla_{\bG} u(z,\sigma)| = c |z|$ on $\partial\Omega$.
Combining \eqref{eq:EL1} and \eqref{eq:EL2}, a critical point of $F_{\bG}$
in the pair $(u, \Omega)$ subject to the boundary conditions $u = 1$ on $\partial B_r(0)$
and $u = 0$ on $\partial\Omega$ is a solution of the overdetermined problem
\begin{equation}
\label{eq:overdetermined}
\begin{cases}
\Delta_{\bG,p} u = 0 & \text{in } \Omega \setminus \overline{B_r(0)}, \\
u = 1 & \text{on } \partial B_r(0), \\
u = 0 & \text{on } \partial\Omega, \\
|\nabla_{\bG} u(z,\sigma)| = c |z| & \text{on } \partial\Omega.
\end{cases}
\end{equation}

\section{The sub-Riemannian Bernoulli problem}
\label{sc:BPmeasure}
Following \cite{GV23} we recall the definitions of partially symmetric domains and of Heisenberg-type $\bG$ that satisfy \emph{Property (H)}.

\begin{definition}
\label{df:partialsy}
Let $\Omega \subset \bG$ be a set in a group of Heisenberg type $\bG$. We say that an open set $\Omega \subset \bG$ has \emph{partial
symmetry} if in the space $\rr \times \rr^k$, with coordinates $(\xi, \sigma)$, there exists an open set $\Omega^{\star}$ such that $(\xi,\sigma) \in \Omega^{\star} \iff (-\xi, \sigma) \in \Omega^{\star} $ and
\[
\Omega=\left\{(z,\sigma) \in \bG \ : \ \left(\frac{|z|^2}{4}, \sigma \right) \in \Omega^{\star}_+\right\},
\]
where we  denote $\Omega^{\star}_+= \{(\xi,\sigma) \in \Omega^{\star} \ : \ \xi \ge 0\}.$
\end{definition}

\begin{definition}\label{D:fparsym}
We say that $u: \bG\to \R$ has \emph{partial symmetry} if, denoting by $(\xi,\sigma)$ the variable point in the space $\R\times \R^k$, there exists $g:\R\times \R^k \to \R$, even in $\xi$, such that
\begin{equation}\label{fgG}
u(z,\sigma) = g\left(\frac{|z|^2}{4},\sigma\right).
\end{equation}
\end{definition}

We indicate with $O(\frak{z}^{\perp})$ the orthogonal group of the vector space $\frak{z}^{\perp}$, where as above $\mathfrak g = \frak{z}^{\perp}  \oplus \ \frak{z}$ is the Lie algebra of $\bG$.  Unlike the standard Laplacian the sub-Laplacian is not invariant with respect to action of $O(\frak{z}^{\perp})$. We denote by $\mathscr S(\frak{z}^{\perp})$ the subgroup of those transformations $S\in O(\frak{z}^{\perp})$ such that for every $\sigma\in \mathfrak{z}$
\begin{equation}\label{commell}
SJ_{\sigma} = J_{\sigma} S.
\end{equation}
We will refer to  $\mathscr S(\frak{z}^{\perp})$ as the $J-$\emph{invariant} orthogonal group of $\frak{z}^{\perp}$ .

\begin{definition}\label{D:trans}
We say that a group of Heisenberg type $\bG$ has the \emph{Property (H)} if the subgroup $\mathscr S(\frak{z}^{\perp})\subset O(\frak{z}^{\perp})$ acts transitively on the unit sphere in $\frak{z}^{\perp}$.
\end{definition}

We shall consider frequently  groups of Heisenberg-type $\bG$ that satisfy \emph{Property (H)}, cf. \cite[Definition 1.4]{GV23}. In \cite[Theorem 5.1]{GV23} it was shown that a group of Heisenberg type satisfies property (H) if and only if it is isomorphic to a (real) Heisenberg group, or a quaternionic Heisenberg group, or a real version of a complexified Heisenberg group.

Invoking a uniqueness argument, the transitivity implies that if the equation and  the domain are invariant under the natural action of $\mathscr S(\frak{z}^{\perp})$ then the solution will be partially symmetric as well. In that sense, the assumption of Property (H) can be dispensed if one is willing to assume instead that the solution of the considered problems has partial symmetry.
Weak solutions of the horizontal $p$-Laplacian with symmetries of this type are studied in \cite{G09grad}: an a priori bound for the right-invariant horizontal gradient gives the $C^{1,\alpha}_{loc}$ regularity of weak solutions with special symmetries; in $\mathbb H^1$ this includes all weak solutions of the form $u(z,t)=u(|z|,t)$, $p\ge2$.

 Recall that  in the setting of groups of Heisenberg type $\bG$ we assume that $p$ is an even integer.  For a function $u$ with partial symmetry defined on an  open set $\Omega \subset \bG$ which has partial
symmetry, cf. \eqref{fgG}, it will be convenient to define  $N=\frac{m+p}{2}\in\nn$, $n=N+k\in\nn$,  and fix unit vectors $e\in\frak z^\perp$, $e'\in\rr^N$.   We let
\begin{equation}\label{e:Omega**, Pi and Psi}
\begin{gathered}
\Omega^{\star \star}= \{ (w,\sigma) \in \rr^{N} \times \rr^k \ : \ (|w|,\sigma) \in \Omega^{\star}_+\},\qquad G(w,\sigma)=g(|w|,\sigma),\\
\Pi(z,\sigma)=(\frac{|z|^2}{4}e',\sigma), \quad \Psi(w,\sigma)=(2\sqrt{|w|}\,e,\sigma).
\end{gathered}
\end{equation}
 By Definition \ref{df:partialsy}, $(z,\sigma)\in\Omega$ if and only if $\Pi(z,\sigma)\in\Omega^{\star\star}$, and $(w,\sigma)\in\Omega^{\star\star}$ if and only if $\Psi(w,\sigma)\in\Omega$.
Hence $\Omega^{\star\star}$ is a bounded open set, $\Psi(\partial\Omega^{\star\star})\subset\partial\Omega$. If $\Omega$ is connected, then $\Omega^{\star\star}$ is connected: it is the union of the $SO(N)$-orbits, in the variable $w$, of the points of the connected set $\Pi(\Omega)$, and each orbit is connected.
Moreover $|\Pi(z,\sigma)|=\frac{|(z,\sigma)|^2_{\bG}}{4}$, $|\Psi(x)|_{\bG}^2=4|x|$, $u=G\circ\Pi$ and $G=u\circ\Psi$, hence we have the relations
\begin{equation}\label{e:u, g and G}
u(z,\sigma)=g(\xi,\sigma)=G(w,\sigma), \quad \xi=\frac{|z|^2}{4}=|w|.
\end{equation}

\begin{proposition}\label{pr:G-weak}
Let $\bG$ have Property $(H)$, $p$ be an even integer, and let $\Omega$ be a bounded open set with partial symmetry, with $\overline{B_r(0)}\subset\Omega$. Let $u\in W^{1,p}_{\bG,0}(\Omega)$, with $u=1$ a.e. in $B_r(0)$, be partially symmetric. Set $\bar r=r^2/4$ and $\Omega'=\Omega^{\star\star}\smallsetminus\overline{B^{\star\star}(0,\bar r)}$.
Then
\begin{itemize}
\item[(i)] $G\in W^{1,p}_0(\Omega^{\star\star})$, $G=1$ a.e. in $B^{\star\star}(0,\bar r)$, $|\nabla_{\bG}u|=\frac{|z|}{2}|\nabla G\circ\Pi|$ a.e. in $\Omega$, where $\Pi$ is defined in \eqref{e:Omega**, Pi and Psi}, and
\[
\kappa\int_{\Omega'}|\nabla G|^p\,dw\,d\sigma=\int_{\Omega\smallsetminus\overline{B_r(0)}}|\nabla_{\bG}u|^p\,dx ,\quad \kappa=\frac{2^{m-1}|S^{m-1}|}{|S^{N-1}|};
\]

\item[(ii)] if, in addition, $u$ is weakly $p$-harmonic in $\Omega\smallsetminus\overline{B_r(0)}$, then $G$ satisfies \eqref{eq:weakdirprob}, $G\in C^{1,\alpha}_{loc}(\Omega^{\star\star}\smallsetminus\overline{B^{\star\star}(0,\bar r)})$ for some $\alpha\in(0,1)$, $u\in C^1(\Omega\smallsetminus\overline{B_r(0)})$ in the Euclidean sense, and the identity in (i) holds at every point of $\Omega\smallsetminus\overline{B_r(0)}$.
\end{itemize}
\end{proposition}

\begin{proof} Let us observe that with $\bar r=\frac{r^2}{4}$ we have  $\overline{B^{\star\star}(0,\bar r)}\subset\Omega^{\star\star}$.

First we prove  (i). Let $\phi_j\in C^\infty_0(\Omega)$, $\phi_j\to u$ in $W^{1,p}_{\bG}(\Omega)$, and $u_j(x)=\int_{\mathscr S(\frak z^\perp)}\phi_j(\Phi_S x)\,dS$, where $dS$ is the normalized Haar measure of the compact group $\mathscr S(\frak z^\perp)$ and $\Phi_S$ is as in the proof of Theorem \ref{th:subapprosimatointerno}. Since $\Phi_S$ preserves $\Omega$, $dx$ and $|\nabla_{\bG}\cdot|$, and $u\circ\Phi_S=u$, we have $u_j\in C^\infty_0(\Omega)$ and $\|u_j-u\|_{W^{1,p}_{\bG}(\Omega)}\le\|\phi_j-u\|_{W^{1,p}_{\bG}(\Omega)}\to0$. By Property $(H)$, $u_j(z,\sigma)=h_j(\frac{|z|^2}{4},\sigma)$, with $h_j(\xi,\sigma)=u_j(2\sqrt\xi\,e,\sigma)$, $\xi\ge0$. Since $-I\in\mathscr S(\frak z^\perp)$, $t\mapsto u_j(te,\sigma)$ is smooth and even, 
so $\partial_tu_j(te,\sigma)=t\,a_j(t,\sigma)$ with $a_j$ smooth, and $\partial_\xi h_j(\xi,\sigma)=2a_j(2\sqrt\xi,\sigma)$. Thus $h_j\in C^1([0,\infty)\times\rr^k)$, and $G_j(w,\sigma)=h_j(|w|,\sigma)$ is Lipschitz, with compact support in $\Omega^{\star\star}$;  in particular $G_j\in W^{1,p}_0(\Omega^{\star\star})$  By \eqref{eq:Xj} and \eqref{eq:Jorth}, for $z\ne0$, $\xi=\frac{|z|^2}{4}$ and $|w|=\xi$,
\begin{equation}\label{eq:polgrad}
\nabla_{\bG}u_j(z,\sigma)=\frac12\Bigl(\partial_\xi h_j\,z+\sum_{\ell=1}^k\partial_{\sigma_\ell}h_j\,J_{\varepsilon_\ell}(z)\Bigr),\qquad
\nabla G_j(w,\sigma)=\Bigl(\partial_\xi h_j\,\frac{w}{|w|},\nabla_\sigma h_j\Bigr),
\end{equation}
with $h_j$ evaluated at $(\xi,\sigma)$. In particular, $|\nabla_{\bG}u_j(z,\sigma)|=\sqrt\xi\,|\nabla G_j(w,\sigma)|$ and $\escpr{\nabla_{\bG}u_j,\nabla_{\bG}u_i}(z,\sigma)=\xi\escpr{\nabla G_j,\nabla G_i}(w,\sigma)$. Writing $z=2\sqrt\xi\,\zeta$, $\zeta\in S^{m-1}$, and $w=\xi\,\omega$, $\omega\in S^{N-1}$, we have $dx=2^{m-1}\xi^{\frac{m-2}2}d\xi\,d\zeta\,d\sigma$ and $dw=\xi^{N-1}d\xi\,d\omega$, with $N-1=\frac{m-2}2+\frac p2$. Hence
\[
\int_{\Omega}|\nabla_{\bG}(u_j-u_i)|^p\,dx=\kappa\int_{\Omega^{\star\star}}|\nabla(G_j-G_i)|^p\,dw\,d\sigma,\qquad
\int_{\Omega^{\star\star}}|G_j-G_i|^p\,dw\,d\sigma\le C\int_{\Omega}|u_j-u_i|^p\,dx,
\]
with $C=\kappa^{-1}\sup_{\Omega}(|z|/2)^p$. Thus $\{G_j\}_{j=1}^\infty$ is a Cauchy sequence in $W^{1,p}_0(\Omega^{\star\star})$. The same change of variables shows that $G_j\to G$ in $L^p(\Omega^{\star\star})$, so $G\in W^{1,p}_0(\Omega^{\star\star})$ and $G_j\to G$ in $W^{1,p}$. Moreover $(\partial_\xi h_j,\nabla_\sigma h_j)$ converges in $L^p(\xi^{N-1}d\xi\,d\sigma)$ to some $(A,B)$, and, passing to the limit in \eqref{eq:polgrad}, $\nabla_{\bG}u=\frac12(A\,z+\sum_\ell B_\ell J_{\varepsilon_\ell}z)$ and $\nabla G=(A\frac{w}{|w|},B)$ a.e. This gives the identity in (i); the integral identity follows by the same change of variables. Finally, $G=1$ a.e. in $B^{\star\star}(0,\bar r)$ because $u=1$ a.e. in $B_r(0)$.

We now prove (ii). Let $\varphi\in C_0^\infty(\Omega')$ be radial in $w$. Then $\psi=\varphi\circ\Pi$ is smooth, partially symmetric and compactly supported in $\Omega\smallsetminus\overline{B_r(0)}$. Indeed, a smooth function radial in $w$ has a smooth even signed radial profile. By \eqref{eq:polgrad} for $\psi$ and the representation of $\nabla_{\bG}u$ above, $\escpr{\nabla_{\bG}u,\nabla_{\bG}\psi}(z,\sigma)=\xi\escpr{\nabla G,\nabla\varphi}(w,\sigma)$ and $|\nabla_{\bG}u|^{p-2}=\xi^{\frac{p-2}{2}}|\nabla G|^{p-2}$ a.e., $|w|=\xi$. Since $\xi^{\frac{m-2}{2}}\xi^{\frac{p-2}{2}}\xi=\xi^{N-1}$,
\[
0=\int_{\Omega\smallsetminus\overline{B_r(0)}}|\nabla_{\bG}u|^{p-2}\langle\nabla_{\bG}u,\nabla_{\bG}\psi\rangle\,dx
=\kappa\int_{\Omega'}|\nabla G|^{p-2}\langle\nabla G,\nabla\varphi\rangle\,dw\,d\sigma .
\]

For a general $\varphi\in C^\infty_0(\Omega')$ set $\bar\varphi(w,\sigma)=\int_{SO(N)}\varphi(Rw,\sigma)\,dR$, which is smooth, compactly supported in $\Omega'$ and radial in $w$.
Write $F=|\nabla G|^{p-2}\nabla G$ and $I(\phi)=\int_{\Omega'}\langle F,\nabla\phi\rangle\,dw\,d\sigma$.
By (i), $|F|=|\nabla G|^{p-1}\in L^{p'}(\Omega')\subset L^1(\Omega')$, since $\Omega'$ is bounded.
For $R\in SO(N)$ let $\rho_R(w,\sigma)=(Rw,\sigma)$, an orthogonal map of $\R^{N+k}$ with $\det\rho_R=1$, leaving $\Omega'$ invariant.

First notice that $I(\varphi\circ\rho_R)=I(\varphi)$ for every $R$.
Since $G$ is radial in $w$, $G\circ\rho_R=G$, hence $\nabla G(\rho_Rx)=\rho_R\nabla G(x)$ a.e., i.e.\ $F(\rho_Rx)=\rho_RF(x)$.
Also $\nabla(\varphi\circ\rho_R)(x)=\rho_R^{T}\nabla\varphi(\rho_Rx)$. Therefore
\[
\langle F(x),\nabla(\varphi\circ\rho_R)(x)\rangle
=\langle\rho_RF(x),\nabla\varphi(\rho_Rx)\rangle
=\langle F(\rho_Rx),\nabla\varphi(\rho_Rx)\rangle ,
\]
and the change of variables $x'=\rho_Rx$ gives $I(\varphi\circ\rho_R)=I(\varphi)$.
Then, let $K=\operatorname{supp}\varphi$ and $K'=\{x:\ \rho_Rx\in K\text{ for some }R\in SO(N)\}$. The set $K'$ is compact
(image of the compact set $SO(N)\times K$ under $(R,y)\mapsto\rho_R^{-1}y$) and $K'\subset\Omega'$.
Differentiating under the integral sign (everything is smooth and bounded),
\[
\bar\varphi(x)=\int_{SO(N)}\varphi(\rho_Rx)\,dR,\qquad
\nabla\bar\varphi(x)=\int_{SO(N)}\nabla(\varphi\circ\rho_R)(x)\,dR .
\]
The function $(x,R)\mapsto\langle F(x),\nabla(\varphi\circ\rho_R)(x)\rangle$ is measurable on $\Omega'\times SO(N)$ and
\[
\bigl|\langle F(x),\nabla(\varphi\circ\rho_R)(x)\rangle\bigr|\le\|\nabla\varphi\|_\infty\,|F(x)|\,\mathbf 1_{K'}(x),
\]
which is integrable on $\Omega'\times SO(N)$ (the measure of $SO(N)$ is $1$ and $F\in L^1(K')$).
Hence Fubini's theorem applies and gives
\[
I(\bar\varphi)=\int_{\Omega'}\Bigl\langle F(x),\int_{SO(N)}\nabla(\varphi\circ\rho_R)(x)\,dR\Bigr\rangle dx
=\int_{SO(N)}\underbrace{\int_{\Omega'}\langle F,\nabla(\varphi\circ\rho_R)\rangle\,dx}_{=\,I(\varphi\circ\rho_R)\,=\,I(\varphi)}\,dR
=I(\varphi).
\]
Hence, the function $\bar\varphi$ is smooth, supported in $K'\subset\Omega'$ and radial in $w$ (because $SO(N)$ acts transitively on the spheres $\{|w|=s\}$, $N\ge2$). By the radial case $I(\bar\varphi)=0$, so $I(\varphi)=0$.
Thus $G$ satisfies \eqref{eq:weakdirprob} for $\varphi\in C^\infty_0(\Omega')$, and, since $|\nabla G|^{p-1}\in L^{p'}(\Omega')$, by density for $\varphi\in W^{1,p}_0(\Omega')$.  By the interior regularity of weak solutions of the Euclidean $p$-Laplace equation \cite{DiB83,Lewis,To84}, $G\in C^{1,\alpha}_{loc}(\Omega')$. Since $u=G\circ\Pi$ a.e. and $\Pi$ is smooth with $\Pi(\Omega\smallsetminus\overline{B_r(0)})\subset\Omega'$, $u$ has a representative in $C^1(\Omega\smallsetminus\overline{B_r(0)})$, whose classical derivatives $X_ju$ coincide a.e. with the distributional ones. The identity in (i) then holds everywhere by continuity.

\end{proof}

\begin{theorem}
\label{th:subapprosimatointerno}
Suppose that $\bG$ has the property $(H)$ and let $p \in \nn$ be an even number. Let $\Omega \subset \bG$ be a connected, bounded open set with partial symmetry such that $0 \in \Omega$ and let $r>0$ such that $B_r(0) \Subset \Omega$. Let $u$ be a  weak solution of
\begin{equation}
\label{eq:BerSub}
 \begin{cases}
        \Delta_{ \bG,p} u=0 &  \Omega \smallsetminus  \overline{B_r(0)},\\
        u=0 & \partial \Omega,\\
        u=1 &\overline{ B_r(0)},\\
    \end{cases}
\end{equation}
such that $u>0$ in $\Omega \smallsetminus  \overline{B_r(0)}$ and assume that there exists $c>0$ such that for each $\eps >0$ there exists a neighborhood $O=O(\eps) \supset \partial \Omega$ such that
\begin{equation}
\label{eq:sRocnb}
\bigl| \, |\nabla_{\bG} u (z,\sigma)|-c|z|\, \bigr| \le \eps |z| , \qquad \textrm{and}\;\;u(z,\sigma)<\eps
\end{equation}
for each $(z,\sigma) \in O \cap \Omega$.
Then, there exists $R>r$ such that $\Omega=B_R(0)$ and
\begin{equation}
\label{eq:usol}
u(z,\sigma)=\begin{cases}
\vspace{0.3 cm}
\displaystyle \frac{|(z, \sigma)|_{\bG}^{\frac{p- Q}{p-1}} - R^{\frac{p- Q}{p-1}}}{r^{\frac{p- Q}{p-1}}- R^{\frac{p- Q}{p-1}}} & p \neq Q,\\
\displaystyle \frac{\ln(\frac{|(z,\sigma)|_{\bG}}{R})}{\ln(\frac{r}{R}) } & p=Q.
\end{cases}
\end{equation}
\end{theorem}

\begin{proof}
We first show that $u$ is partially symmetric. Let $S\in\mathscr S(\frak z^\perp)$ and $\Phi_S(z,\sigma)=(Sz,\sigma)$. By \eqref{eq:J} and \eqref{commell}, for $z,z'\in\frak z^\perp$ and $\sigma\in\frak z$,
\[
\escpr{\sigma,[Sz,Sz']}=\escpr{J_\sigma Sz,Sz'}=\escpr{SJ_\sigma z,Sz'}=\escpr{J_\sigma z,z'}=\escpr{\sigma,[z,z']},
\]
so that $[Sz,Sz']=[z,z']$ and, by the group law, $\Phi_S$ is a group automorphism. Since $\Phi_S$ is linear and orthogonal, it preserves $dx$. Moreover $\Phi_S\bigl(x\circ(te_j,0)\bigr)=\Phi_S(x)\circ(tSe_j,0)$, hence for $f\in C^1$
\[
\nabla_{\bG}(f\circ\Phi_S)=S^{t}\,(\nabla_{\bG}f)\circ\Phi_S ,
\]
and by density the same holds for $f\in W^{1,p}_{\bG}$. Since $|Sz|=|z|$, $\Phi_S$ maps $B_r(0)$, $\Omega$ and $\Omega\smallsetminus\overline{B_r(0)}$ onto themselves, and $C^\infty_0(\Omega)$ onto itself; thus $u\circ\Phi_S\in W^{1,p}_{\bG,0}(\Omega)$ and $u\circ\Phi_S=1$ in $B_r(0)$. For $\varphi\in C^\infty_0(\Omega\smallsetminus\overline{B_r(0)})$, the change of variables $x'=\Phi_S(x)$ and the orthogonality of $S$ give
\[
\int|\nabla_{\bG}(u\circ\Phi_S)|^{p-2}\escpr{\nabla_{\bG}(u\circ\Phi_S),\nabla_{\bG}\varphi}\,dx
=\int|\nabla_{\bG}u|^{p-2}\escpr{\nabla_{\bG}u,\nabla_{\bG}(\varphi\circ\Phi_S^{-1})}\,dx'=0 .
\]
Hence $u\circ\Phi_S$ is a weak solution of \eqref{eq:BerSub}, and the comparison principle \cite[Lemma 2.6]{D95} gives $u\circ\Phi_S=u$; it also gives $0\le u\le1$. Since weak solutions are continuous \cite{CDG96}, $u(Sz,\sigma)=u(z,\sigma)$ for every $(z,\sigma)\in\Omega$ and $S\in\mathscr S(\frak z^\perp)$. If $|z|=|z'|$, Property $(H)$ gives $S$ with $z'=Sz$. Hence $u(z,\sigma)=\hat u(|z|,\sigma)$, and $g(\xi,\sigma)=\hat u(2\sqrt{|\xi|},\sigma)$ is even in $\xi$ and satisfies $u(z,\sigma)=g(\frac{|z|^2}{4},\sigma)$.

Set $\bar r=r^2/4$. Since $\overline{B_r(0)}\subset\Omega$, the identities in \eqref{e:Omega**, Pi and Psi} give $\overline{B^{\star\star}(0,\bar r)}\subset\Omega^{\star\star}$.

By Proposition \ref{pr:G-weak}, $G\in W^{1,p}_0(\Omega^{\star\star})$, $G=1$ a.e. in $\overline{B^{\star\star}(0,\bar r)}$, $G$ satisfies \eqref{eq:weakdirprob}, $G\in C^{1,\alpha}_{loc}(\Omega^{\star\star}\smallsetminus\overline{B^{\star\star}(0,\bar r)})$, $X_1u,\dots,X_mu$ are continuous in $\Omega\smallsetminus\overline{B_r(0)}$, and
\begin{equation}\label{eq:gradcorr}
|\nabla_{\bG}u(z,\sigma)|=\frac{|z|}{2}\,|\nabla G(\Pi(z,\sigma))|,\qquad (z,\sigma)\in\Omega\smallsetminus\overline{B_r(0)}.
\end{equation}
Since $u>0$, $G>0$ in $\Omega^{\star\star}\smallsetminus\overline{B^{\star\star}(0,\bar r)}$.

Let $0<\eps<1$ and $O=O(\eps)$. Since $\partial\Omega$ is compact, there is $\delta>0$ such that every $y\in\Omega$ with $\operatorname{dist}(y,\partial\Omega)<\delta$ lies in $O$. Since $\Psi$ is uniformly continuous on bounded sets, there is $\eta>0$ such that $|\Psi(x)-\Psi(x')|<\delta$ whenever $x,x'\in\overline{\Omega^{\star\star}}$ and $|x-x'|<\eta$. Let $x=(w,\sigma)\in\Omega^{\star\star}$ with $\operatorname{dist}(x,\partial\Omega^{\star\star})<\eta$. Since $\Psi(\partial\Omega^{\star\star})\subset\partial\Omega$, $\Psi(x)=(z,\sigma)\in O\cap\Omega$, with $|z|=2\sqrt{|w|}$. By \eqref{eq:sRocnb} and \eqref{eq:gradcorr},
\[
\frac{|z|}{2}\,\bigl|\,|\nabla G(x)|-2c\,\bigr|\le\eps|z|,\qquad G(x)=u(\Psi(x))<\eps .
\]
Hence $\bigl||\nabla G(x)|-2c\bigr|\le2\eps$ if $w\ne0$, and, by continuity of $\nabla G$, also if $w=0$. Thus \eqref{ovrdtrmnd} holds, with $2r(\eps)=\eta(\eps/2)$. By Theorem \ref{th:euclideointerno} there is $\bar R>\bar r$ such that $\Omega^{\star\star}=B^{\star\star}(0,\bar R)$ and
\[
 G(x)=\begin{cases}
 \displaystyle \frac{|x|^{\frac{p-n}{p-1}}- \bar R^{\frac{p-n}{p-1}} }{\bar r^{\frac{p-n}{p-1}}- \bar R^{\frac{p-n}{p-1}}} & p\ne n,\\[10pt]
 \displaystyle \frac{\ln(\frac{|x|}{\bar R})}{\ln(\frac{\bar r}{\bar R})} & p=n.
 \end{cases}
\]
Since $(z,\sigma)\in\Omega$ if and only if $\frac{|(z,\sigma)|^2_{\bG}}{4}<\bar R$, $\Omega=B_R(0)$ with $R=2\bar R^{\frac12}$. Since $n=\frac{Q+p}{2}$, $p=n$ if and only if $p=Q$, and $\frac{p-n}{p-1}=\frac{p-Q}{2(p-1)}$. Hence, by $u=G\circ\Pi$ and $\bar r=\frac{r^2}{4}$, for $p\ne Q$
\[
u(z,\sigma)=\frac{\bigl(\frac{|(z,\sigma)|_{\bG}^2}{4}\bigr)^{\frac{p-Q}{2(p-1)}}-\bigl(\frac{R^2}{4}\bigr)^{\frac{p-Q}{2(p-1)}}}{\bigl(\frac{r^2}{4}\bigr)^{\frac{p-Q}{2(p-1)}}-\bigl(\frac{R^2}{4}\bigr)^{\frac{p-Q}{2(p-1)}}}
=\frac{|(z,\sigma)|_{\bG}^{\frac{p-Q}{p-1}} - R^{\frac{p-Q}{p-1}}}{r^{\frac{p-Q}{p-1}}-R^{\frac{p-Q}{p-1}}},
\]
and for $p=Q$
\[
u(z,\sigma)=\frac{\ln\bigl(|(z,\sigma)|_{\bG}^2/R^2\bigr)}{\ln(r^2/R^2)}
=\frac{\ln\bigl(|(z,\sigma)|_{\bG}/R\bigr)}{\ln(r/R)},
\]
which is \eqref{eq:usol}.

\end{proof}

\begin{remark}
    Notice that if $\frac{|\nabla_{\bG} u(z,\sigma)|}{|z|}=c$ for each $(z,\sigma) \in \partial \Omega$ and the function $(z,\sigma) \to \frac{|\nabla_{\bG} u(z,\sigma)|}{|z|}$ extends continuously to $\overline\Omega\smallsetminus B_r(0)$, and $u\in C(\overline\Omega\smallsetminus B_r(0))$ with $u=0$ on $\partial\Omega$, then the condition \eqref{eq:sRocnb} follows by continuity.
\end{remark}

\begin{remark}\label{rm:pLap}
The choice $n=\frac{m+p}{2}+k$ in the proof of Theorem~\ref{th:subapprosimatointerno} has the following pointwise counterpart, cf. \cite[Lemma 2.4 and Proposition 2.6]{GV23} and \cite[Propositions 3.2 \& 3.3]{GVduke01}. Let $u(z,\sigma)=g(\xi,\sigma)$, $\xi=\frac{|z|^2}{4}>0$, with $g\in C^2$, and let $G(w,\sigma)=g(|w|,\sigma)$, $(w,\sigma)\in\rr^N\times\rr^k$, $N=\frac{m+p}{2}$. By \eqref{eq:Xj} and \eqref{eq:Jorth}, $|\nabla_{\bG}u|^2=\xi(g_\xi^2+|\nabla_\sigma g|^2)$ and $\Delta_{\bG}u=\xi(g_{\xi\xi}+\frac m{2\xi}g_\xi+\Delta_\sigma g)$, and a direct computation gives
\[
\Delta_{\bG,p}u=\xi^{\frac p2}\Bigl\{|\nabla g|^{p-2}\Bigl(g_{\xi\xi}+\frac{m+p-2}{2\xi}g_\xi+\Delta_\sigma g\Bigr)+(p-2)|\nabla g|^{p-4}\Delta_\infty g\Bigr\},
\]
where $|\nabla g|^2=g_\xi^2+|\nabla_\sigma g|^2$, $\Delta_\infty g=\frac12\langle\nabla|\nabla g|^2,\nabla g\rangle$, and the last term is absent for $p=2$. Since $\frac{m+p-2}{2}=N-1$, the expression in braces is the Euclidean $p$-Laplacian of $G$ in $\rr^n$ at $(w,\sigma)$, $|w|=\xi$. For even $p$ it is a polynomial in the derivatives of $g$, hence meaningful at critical points of $g$. The proofs do not use this identity, but its weak form, Proposition~\ref{pr:G-weak}.
\end{remark}

\section{Euclidean result}

In this section we prove the Euclidean rigidity result used in the proof of Theorem~\ref{th:subapprosimatointerno}. The function $G$ is only assumed to be a weak solution, and no regularity of $\partial\Omega^{\star\star}$ is assumed. The proof compares $G$, on its level sets near $\partial\Omega^{\star\star}$, with the radial $p$-capacitary potentials of the largest inscribed and the smallest circumscribed balls centered at the origin.

\begin{theorem}
\label{th:euclideointerno}
Let $1<p<+\infty$ and $n \in \nn$ with $n \ge 2$. Let $\Omega^{\star \star} \subset \rr^n$ be  a connected, bounded open set with $\overline{B^{\star \star}(0,\bar r)}\subset\Omega^{\star \star}$ and suppose that $G \in W_0^{1,p}(\Omega^{\star \star} )$ verifies
    \begin{equation}
    \label{eq:weakdirprob}
        \int_{ \Omega^{\star \star} \smallsetminus  \overline{B^{\star \star}(0,\bar r)} } |\nabla G|^{p-2}\escpr{ \nabla G, \nabla \varphi} dw \,d\sigma =0 ,
    \end{equation}
  for each $\varphi \in W_0^{1,p}(\Omega^{\star \star} \smallsetminus  \overline{B^{\star \star}(0,\bar r)})$ and $G=1$ a.e. in $\overline{B^{\star \star}(0,\bar r)}$.
  Assume that $G>0$ in $\Omega^{\star \star} \smallsetminus  \overline{B^{\star \star}(0,\bar r)}$ and that there exists $c>0$ such that for every $0<\eps <1$ we have
  \begin{equation}\label{ovrdtrmnd} \big| |\nabla G(w,\sigma)| -2c\big| \le \eps, \qquad \textrm{and}\;\; G(w,\sigma)<\eps
  \end{equation}
  for each $(w,\sigma)\in\Omega^{\star\star}$ with $\operatorname{dist}((w,\sigma),\partial\Omega^{\star\star})<2r(\eps)$, for some $r(\eps)>0$.
Then there exists $\bar R$ such that $\Omega^{\star \star}$ is a Euclidean ball $B^{\star \star}(0, \bar R)$ and
 \[
 G(x)=\begin{cases} \vspace{0.3cm}
 \displaystyle \frac{|x|^{\frac{p-n}{p-1}}- \bar R^{\frac{p-n}{p-1}} }{\bar r^{\frac{p-n}{p-1}}- \bar R^{\frac{p-n}{p-1}}} & p\ne n,\\

 \displaystyle \frac{\ln(\frac{|x|}{\bar R})}{\ln(\frac{\bar r}{\bar R})} & p=n.
 \end{cases}
 \]
\end{theorem}

\begin{proof}
Let $R_1=\operatorname{dist}(0,\partial\Omega^{\star\star})$, $R_2=\max_{\partial\Omega^{\star\star}}|x|$, and $x_i\in\partial\Omega^{\star\star}$ with $|x_i|=R_i$. Then $B^{\star\star}(0,R_1)\subset\Omega^{\star\star}\subset B^{\star\star}(0,R_2)$, $R_1>\bar r$, and $R_1\le R_2$, with equality if and only if $\Omega^{\star\star}$ is a ball. We prove $R_2\le R_1$.

Let $A=\Omega^{\star\star}\smallsetminus\overline{B^{\star\star}(0,\bar r)}$, and let $G=1$ on $\overline{B^{\star\star}(0,\bar r)}$. By \cite{DiB83,Lewis,To84}, $G\in C^{1,\alpha}_{loc}(A)$. For $\rho>\bar r$ let
\[
v_\rho(x)=\frac{\Lambda(|x|)-\Lambda(\rho)}{\Lambda(\bar r)-\Lambda(\rho)},\qquad \Lambda(s)=s^{\frac{p-n}{p-1}}\ \ (p\ne n),\qquad \Lambda(s)=\ln s\ \ (p=n),
\]
the radial $p$-harmonic function with $v_\rho=1$ on $\partial B^{\star\star}(0,\bar r)$ and $v_\rho=0$ on $\partial B^{\star\star}(0,\rho)$. On $\partial B^{\star\star}(0,\rho)$, $|\nabla v_\rho|=\kappa_{n,p}f(\rho)$, where $f$ is given by
	    \begin{equation}
        \label{def:f}
	    f(t)=
	    \begin{cases}
	    \displaystyle \frac{t^{\frac{1-n}{p-1}}}{\bar{r}^{\frac{p-n}{p-1}}- t^{\frac{p-n}{p-1}}} & p<n\\
	    \displaystyle \frac{t^{\frac{1-n}{p-1}}}{t^{\frac{p-n}{p-1}}- \bar{r}^{\frac{p-n}{p-1}}} & p>n\\
	    \displaystyle \frac{1}{t \ln (\frac{t}{\bar r})} & p=n
	    \end{cases}
	    \end{equation}
 and $\kappa_{n,p}=\frac{|n-p|}{p-1}$ for $p\ne n$ and $\kappa_{n,p}=1$ for $p=n$.

Since $0\le(G-1)^+\le G$, $(G-1)^+\in W^{1,p}_0(A)$, and testing \eqref{eq:weakdirprob} with it gives $G\le1$. Let $\tilde v=\min\{v_{R_1},1\}$ in $B^{\star\star}(0,R_1)$. Then $0\le(\tilde v-G)^+\le\tilde v$ and $(\tilde v-G)^+=0$ in $B^{\star\star}(0,\bar r)$, so $(\tilde v-G)^+\in W^{1,p}_0(B^{\star\star}(0,R_1)\smallsetminus\overline{B^{\star\star}(0,\bar r)})$. Testing the equations of $G$ and $v_{R_1}$ with it, the strict monotonicity of $\zeta\mapsto|\zeta|^{p-2}\zeta$ gives $G\ge v_{R_1}$ in $B^{\star\star}(0,R_1)\smallsetminus\overline{B^{\star\star}(0,\bar r)}$. In particular $G$ is continuous in $\Omega^{\star\star}$.

Fix $0<\eps<\min\{1,2c\}$, and let $r(\eps)>0$ be the number given by the hypothesis \eqref{ovrdtrmnd} corresponding to such $\eps$. Decreasing $r(\eps)$ if necessary, we may assume that $r(\eps)<(R_1-\bar r)/2$; in particular, the compact set $K$ defined below is nonempty. Let $\rho_0=\frac{\bar r+R_1}{2}$ and
\[
K=\{x\in\Omega^{\star\star}:\operatorname{dist}(x,\partial\Omega^{\star\star})\ge r(\eps),\ |x|\ge\rho_0\},\qquad m(\eps)=\min\Bigl\{\frac{\Lambda(\rho_0)-\Lambda(R_1)}{\Lambda(\bar r)-\Lambda(R_1)},\ \min_K G\Bigr\}>0 .
\]
$K$ is a compact subset of $A$. If $x\in A$ and $G(x)<m(\eps)$, then $|x|>\rho_0$, because $G\ge v_{R_1}$, and $x\notin K$. Hence $\operatorname{dist}(x,\partial\Omega^{\star\star})<r(\eps)$, and by \eqref{ovrdtrmnd}
\[
0<2c-\eps\le|\nabla G(x)|\le2c+\eps .
\]
Let $a_j\downarrow0$, $0<a_j<\min\{m(\eps),\eps\}$, and $\Omega^{\star\star}_j=\{x\in\Omega^{\star\star}:G(x)>a_j\}$, an open set containing $\overline{B^{\star\star}(0,\rho_0)}$. By \eqref{ovrdtrmnd} with $a_j$ in place of $\eps$, and $r(a_j)>0$ the corresponding number, $G<a_j$ where $\operatorname{dist}(x,\partial\Omega^{\star\star})<2r(a_j)$. Hence $\overline{\Omega^{\star\star}_j}\subset\Omega^{\star\star}$, $\partial\Omega^{\star\star}_j\subset\{x\in A:G(x)=a_j\}$, and $2c-\eps\le|\nabla G|\le2c+\eps$ on $\partial\Omega^{\star\star}_j$. Since $G>0$ in $\Omega^{\star\star}$, every compact subset of $\Omega^{\star\star}$ lies in $\Omega^{\star\star}_j$ for $j$ large. Let $R^j_1=\operatorname{dist}(0,\partial\Omega^{\star\star}_j)$, $R^j_2=\max_{\partial\Omega^{\star\star}_j}|x|$, and $x^j_i\in\partial\Omega^{\star\star}_j$ with $|x^j_i|=R^j_i$. Then $B^{\star\star}(0,R^j_1)\subset\Omega^{\star\star}_j\subset B^{\star\star}(0,R^j_2)$, $\rho_0\le R^j_1\le R_1$, $R^j_2\le R_2$, and $R^j_i\to R_i$ as $j\to\infty$.

For $i=1,2$ let $v^j_i=a_j+(1-a_j)v_{R^j_i}$. Then $v^j_i=1$ on $\partial B^{\star\star}(0,\bar r)$, $v^j_i=a_j$ on $\partial B^{\star\star}(0,R^j_i)$, and $|\nabla v^j_i|=(1-a_j)\kappa_{n,p}f(R^j_i)$ on $\partial B^{\star\star}(0,R^j_i)$. In $U_1=B^{\star\star}(0,R^j_1)\smallsetminus\overline{B^{\star\star}(0,\bar r)}$, $G$ is continuous up to $\partial U_1$, $G\ge a_j=v^j_1$ on $\partial B^{\star\star}(0,R^j_1)$ and $G=v^j_1=1$ on $\partial B^{\star\star}(0,\bar r)$. For $\delta>0$, $(v^j_1-G-\delta)^+$ has compact support in $U_1$, and testing the two equations with it gives $G\ge v^j_1-\delta$. Hence $G\ge v^j_1$ in $U_1$. In the same way, in $U_2=\Omega^{\star\star}_j\smallsetminus\overline{B^{\star\star}(0,\bar r)}$ we have $G=a_j\le v^j_2$ on $\partial\Omega^{\star\star}_j$ and $G=v^j_2=1$ on $\partial B^{\star\star}(0,\bar r)$, so $G\le v^j_2$ in $U_2$.

At $x^j_1$, $G(x^j_1)=a_j=v^j_1(x^j_1)$ and $tx^j_1\in U_1$ for $t<1$ close to $1$. Since $G$ is $C^1$ near $x^j_1$,
\[
|\nabla v^j_1(x^j_1)|=\lim_{t\to1^-}\frac{v^j_1(tx^j_1)-a_j}{(1-t)|x^j_1|}\le\lim_{t\to1^-}\frac{G(tx^j_1)-a_j}{(1-t)|x^j_1|}\le|\nabla G(x^j_1)|\le2c+\eps .
\]
At $x^j_2$, let $\nu=\nabla G(x^j_2)/|\nabla G(x^j_2)|$. For $t>0$ small, $x^j_2+t\nu\in U_2$, and
\[
|\nabla v^j_2(x^j_2)|\ge\lim_{t\to0^+}\frac{v^j_2(x^j_2+t\nu)-a_j}{t}\ge\lim_{t\to0^+}\frac{G(x^j_2+t\nu)-a_j}{t}=|\nabla G(x^j_2)|\ge2c-\eps .
\]
That is,
\begin{equation}\label{Kaza1}
(1-a_j)\kappa_{n,p}f(R^j_1)\le2c+\eps,
\end{equation}
\begin{equation}\label{Kaza2}
(1-a_j)\kappa_{n,p}f(R^j_2)\ge2c-\eps .
\end{equation}
We set $v_i=v_{R_i}$. Combining \eqref{Kaza1} and \eqref{Kaza2} we get
     $$
     \begin{array}{cc}
    (1-a_j)\frac{n-p}{p-1} \frac{(R^j_1)^{\frac{1-n}{p-1}}}{\bar{r}^{\frac{p-n}{p-1}}- (R^j_1)^{\frac{p-n}{p-1}}}   \le  (1-a_j)\frac{n-p}{p-1} \frac{(R^j_2)^{\frac{1-n}{p-1}}}{\bar{r}^{\frac{p-n}{p-1}}- (R^j_2)^{\frac{p-n}{p-1}}}+2\eps
   ,\quad \quad
&\textrm{for}\;  p< n,
 \\
    (1-a_j) \frac{p-n}{p-1} \frac{(R^j_1)^{\frac{1-n}{p-1}}}{(R^j_1)^{\frac{p-n}{p-1}} - \bar{r}^{\frac{p-n}{p-1}} }  \le  (1-a_j)\frac{p-n}{p-1} \frac{(R^j_2)^{\frac{1-n}{p-1}}}{(R^j_2)^{\frac{p-n}{p-1}} - \bar{r}^{\frac{p-n}{p-1}} } +2\eps
  ,\quad \quad
&\textrm{for}\; p>n,
    \\
       \frac{(1-a_j)}{R^j_1 \ln(\frac{R^j_1}{\bar r})} \le  \frac{(1-a_j)}{R^j_2 \ln(\frac{R^j_2}{\bar r})} +2\eps
       ,\quad \quad
&\textrm{for}\; p= n.
    \end{array}
    $$
    Passing first to the limit as $j\to+\infty$ and then letting $\eps \to 0$, we deduce
         \begin{equation}
         \label{Kaza3}
     \begin{array}{cc}
    \frac{(R_1)^{\frac{1-n}{p-1}}}{\bar{r}^{\frac{p-n}{p-1}}- (R_1)^{\frac{p-n}{p-1}}}   \le  \frac{(R_2)^{\frac{1-n}{p-1}}}{\bar{r}^{\frac{p-n}{p-1}}- (R_2)^{\frac{p-n}{p-1}}}
       ,\quad \quad
&\textrm{for}\;
    p< n,
     \\
     \frac{(R_1)^{\frac{1-n}{p-1}}}{(R_1)^{\frac{p-n}{p-1}} - \bar{r}^{\frac{p-n}{p-1}} }  \le   \frac{(R_2)^{\frac{1-n}{p-1}}}{(R_2)^{\frac{p-n}{p-1}} - \bar{r}^{\frac{p-n}{p-1}} }
      ,\quad \quad
&\textrm{for}\; p>n,
    \\
      \frac{1}{R_1 \ln(\frac{R_1}{\bar r})} \le  \frac{1}{R_2 \ln(\frac{R_2}{\bar r})}
     ,\quad \quad
&\textrm{for}\; p=n.
    \end{array}
    \end{equation}
Since the function $f:(\bar r, + \infty) \to \rr$  defined in \eqref{def:f}
   is decreasing in $(\bar r, + \infty)$, by \eqref{Kaza3} we get $R_2\le R_1$.
    Hence $R_2=R_1=\bar R$ and $\Omega^{\star \star}=B^{\star \star}(0,\bar R)$. Notice that $\bar R>0$ depends only on $c$ and $\bar r$ and it is the unique solution of the equations:
    \begin{equation}
    \label{eq:Rcr}
      \varphi(\rho)=\rho^{\frac{p-n}{p-1}} + \frac{n-p}{p-1}\frac{\rho ^{\frac{1-n}{p-1}}}{2c} - \bar{r}^{\frac{p-n}{p-1}}=0
    \end{equation}
    when $p<n$,
        \begin{equation}
    \label{eq:Rcr2}
      \psi(\rho)=\rho^{\frac{p-n}{p-1}} - \frac{p-n}{p-1}\frac{\rho ^{\frac{1-n}{p-1}}}{2c} - \bar{r}^{\frac{p-n}{p-1}}=0
    \end{equation}
    when $p>n$,
    \begin{equation}
    \label{eq:Rcr3}
    \theta(\rho)= \rho \ln(\frac{\rho}{\bar r})-\frac{1}{2c}=0
    \end{equation}
    when $p=n$.
 Indeed, since $\varphi'(\rho)=\frac{p-n}{p-1} \rho^{\frac{1-n}{p-1}}+ \frac{n-p}{p-1}\frac{1-n}{p-1} \frac{\rho ^{\frac{2-n-p}{p-1}}}{2c}<0$ for all $\rho>0$ (notice $ p<n$), and since $\varphi$ is continuous, $\varphi(\bar r)>0$, and $\lim_{\rho\to +\infty} \varphi(\rho)<0$, there exists a unique solution $\bar R>\bar r>0$ of the equation \eqref{eq:Rcr}. Since $\psi'(\rho)=\frac{p-n}{p-1} \rho^{\frac{1-n}{p-1}}+ \frac{p-n}{p-1} \frac{n-1}{p-1} \frac{\rho ^{\frac{2-n-p}{p-1}}}{2c}>0$  for all $\rho>0$ (notice $ p>n$), $\psi(\bar r)<0$ and $\lim_{\rho\to +\infty} \psi(\rho)=+\infty$ there exists a unique solution $\bar R>\bar r>0$ of the equation \eqref{eq:Rcr2}. Since $\theta (\bar r)<0$, $\theta'(\rho)=\ln(\frac{\rho}{\bar r}) +1>0$  and $\lim_{\rho\to +\infty} \theta(\rho)=+\infty$ we get there exists a unique solution $\bar R>\bar r>0$ of the equation \eqref{eq:Rcr3}.
    Moreover, by the uniqueness of the solution of the Dirichlet problem \eqref{eq:weakdirprob} in $ B^{\star \star}(0,\bar R) \smallsetminus B^{\star \star}(0,\bar r)$ we obtain that
    \[
    G(x)=v_1(x)=v_2(x)=\begin{cases} \vspace{0.3cm}
 \displaystyle \frac{|x|^{\frac{p-n}{p-1}}- \bar R^{\frac{p-n}{p-1}} }{\bar r^{\frac{p-n}{p-1}}- \bar R^{\frac{p-n}{p-1}}}, & p\ne n,\\

 \displaystyle \frac{\ln(\frac{|x|}{\bar R})}{\ln(\frac{\bar r}{\bar R})}, & p=n.
 \end{cases}
    \]
\end{proof}

\begin{remark}
    A boundary condition similar to \eqref{ovrdtrmnd} was considered in \cite[Theorem 1]{GL89}. Note that their assumption holds almost everywhere, as it relies on an integral technique; here, since we adopt a pointwise technique, we correspondingly require a pointwise assumption.
\end{remark}

\begin{remark}\label{rm:interiorBP}
The technique of inscribed and circumscribed balls also applies to the interior Bernoulli problem, provided $n\ge3$ and $\Omega^{\star\star}$ contains a priori the ball $B^{\star\star}\bigl(0,\bar R(n-1)^{-\frac{1}{n-2}}\bigr)$. For simplicity let $p=2$. Let $\Omega^{\star\star}$ be a bounded open set with $0\in\Omega^{\star\star}$ and $\overline{\Omega^{\star\star}}\subset B^{\star\star}(0,\bar R)$, and let $G$ solve
\begin{equation*}
 \begin{cases}
        \Delta G=0 &   \text{in } {B^{\star \star}(0,\bar R)} \smallsetminus \overline{\Omega^{\star \star}} ,  \\
        G=1 & \text{on } \partial \Omega^{\star \star},\\
        G=0 & \text{on } \partial {B^{\star \star}(0,\bar R)},
    \end{cases}
\end{equation*}
with the condition \eqref{ovrdtrmnd}, $1-G$ in place of $G$, in an exterior neighbourhood of $\partial\Omega^{\star\star}$. Let $R_1$ be the smallest radius with $\Omega^{\star\star}\subset B^{\star\star}(0,R_1)$, $R_2$ the largest radius with $B^{\star\star}(0,R_2)\subset\Omega^{\star\star}$, and $x_i\in\partial\Omega^{\star\star}\cap\partial B^{\star\star}(0,R_i)$. Then $0<R_2\le R_1<\bar R$, with equality if and only if $\Omega^{\star\star}$ is a ball. Let
\[
v_i(x)=\frac{|x|^{2-n}-\bar R^{2-n}}{R_i^{2-n}-\bar R^{2-n}},\qquad f(t)=\frac{t^{1-n}}{t^{2-n}-\bar R^{2-n}},\quad 0<t<\bar R .
\]
Then $v_i$ is harmonic, $v_i=1$ on $\partial B^{\star\star}(0,R_i)$, $v_i=0$ on $\partial B^{\star\star}(0,\bar R)$, and $|\nabla v_i|=(n-2)f(R_i)$ on $\partial B^{\star\star}(0,R_i)$. By comparison, $G\le v_1$ in $B^{\star\star}(0,\bar R)\smallsetminus\overline{B^{\star\star}(0,R_1)}$ and $G\ge v_2$ in $B^{\star\star}(0,\bar R)\smallsetminus\overline{\Omega^{\star\star}}$. Since $G(x_i)=v_i(x_i)=1$, this gives, after the approximation by level sets of the proof of Theorem \ref{th:euclideointerno},
\[
|\nabla v_1(x_1)|\le|\nabla G(x_1)|=2c,\qquad |\nabla v_2(x_2)|\ge|\nabla G(x_2)|=2c,
\]
hence $f(R_1)\le\frac{2c}{n-2}\le f(R_2)$. Now
\[
f'(t)=\frac{t^{-n}\bigl((n-1)\bar R^{2-n}-t^{2-n}\bigr)}{\bigl(t^{2-n}-\bar R^{2-n}\bigr)^2},
\]
so $f$ is decreasing in $\bigl(0,\bar R(n-1)^{-\frac1{n-2}}\bigr]$ and increasing in $\bigl[\bar R(n-1)^{-\frac1{n-2}},\bar R\bigr)$. If $R_2\ge\bar R(n-1)^{-\frac1{n-2}}$, then $f(R_1)\le f(R_2)$ and $R_2\le R_1$ give $R_1\le R_2$. Hence $R_1=R_2$, and $\Omega^{\star\star}$ is a ball.
\end{remark}

\section{A sharp Euclidean result with  applications to the sub-Riemannian setting}
\label{sc:Bernoullismooth}
The next theorem gives a version of the exterior Bernoulli problem in which no
smoothness of $\partial\Omega^{\star\star}$ is assumed; compare
\cite{ABS21,HP,HPP}. If $\partial\Omega^{\star\star}$ is of class $C^1$ and $G\in C^1(\overline{\Omega^{\star\star}}\smallsetminus B^{\star\star}(0,\bar r))$ with $\nabla G\ne0$ on $\partial\Omega^{\star\star}=\{G=0\}$, then the outer unit normal is $\nu=-\nabla G/|\nabla G|$, since $G>0$ in $\Omega^{\star\star}$. Therefore the overdetermined Neumann boundary condition is
    \begin{equation}
    \label{eq:Neumcond}
      \frac{\partial G(x)}{\partial \nu}=\escpr{\nabla G(x),\nu(x)}=-|\nabla G(x)|=-c,
    \end{equation}
    for $x \in \partial \Omega^{\star \star} $.
    For a boundary point $x\in\partial\Omega^{\star\star}$ at which the interior limit below is meaningful, we formulate the metric boundary-slope condition as
    \begin{equation}
    \label{eq:limsupgradient}
    \limsup_{\substack{y\to x\\y\in\Omega^{\star\star}}}
    \frac{G(y)}{|x-y|}=c.
    \end{equation}
Here and below the limit is taken from inside the domain, using the continuous
interior representative of $G$. For a merely Sobolev function, the condition
$G=0$ on $\partial\Omega^{\star\star}$ is understood in the trace sense and
does not by itself define a pointwise boundary slope. 
In the  next two theorems we shall consider the smallest and largest closed balls centered at the origin that define an annulus containing a given bounded closed set $\Gamma$. Let $V$ be a finite dimensional vector space and $\|\cdot\|: V\to \rr$ a norm. Let
\begin{equation}\label{e:circimscribed and inscribed balls}
\operatorname{Max}(\Gamma)=\operatorname {argmax}_{\Gamma} \|x\|, \quad\text{resp., }\   \operatorname{Min}(\Gamma)=\operatorname {argmin}_{\Gamma} \|x\|,
\end{equation}
be the points of contact of the outer, resp. inner, ball and the set $\Gamma$.

\begin{theorem}
\label{th:eucsharp}
Let $1<p<+\infty$ and $n\ge 2$. Let $\Omega^{\star \star} \subset \rr^n$ be  a connected, bounded open set with $\overline{B^{\star \star}(0,\bar r )}\subset\Omega^{\star \star} $. Suppose that $G \in W_0^{1,p}(\Omega^{\star \star} )$ verifies
    \begin{equation}
        \int_{ \Omega^{\star \star} \smallsetminus  \overline{B^{\star \star}(0,\bar r )} } |\nabla G|^{p-2}\escpr{ \nabla G, \nabla \varphi} \, dx =0 ,\quad \forall\,\varphi \in W_0^{1,p}\left(\Omega^{\star \star} \smallsetminus  \overline{B^{\star \star}(0,\bar r )}\right),
    \end{equation}
 and $G=1$ a.e. in $\overline{B^{\star \star}(0,\bar r )}$.
Assume that for some $c>0$ there exist points $x_1 \in \text{Min}({\partial\Omega}^{**}) $ and $x_2 \in \text{Max}({\partial\Omega}^{**}) $ such that
\begin{equation}
\limsup_{\substack{x\to x_1\\x\in\Omega^{\star\star}}}
\frac{G(x)}{|x-x_1|}\le c
\quad\text{and}\quad
\limsup_{\substack{x\to x_2\\x\in\Omega^{\star\star}}}
\frac{G(x)}{|x-x_2|}\ge c,
\end{equation}
where
\[
\operatorname{Min}({\partial\Omega}^{**})=\operatorname{argmin}_{\partial\Omega^{**}} |x|, \qquad
\operatorname{Max}({\partial\Omega}^{**})= \operatorname{argmax}_{\partial\Omega^{**}} |x|.
\]
 Then there exists $R>\bar r$ such that  $\Omega^{\star \star}$ is a Euclidean ball $B^{\star \star}(0,R)$ and
 \[
 G(x)=\begin{cases} \vspace{0.3cm}
\displaystyle \frac{|x|^{\frac{p-n}{p-1}}- R^{\frac{p-n}{p-1}} }{\bar r^{\frac{p-n}{p-1}}- R^{\frac{p-n}{p-1}}} & p\ne n,\\
\displaystyle \frac{\ln(\frac{|x|}{R})}{\ln(\frac{\bar r}{R})} & p=n.
 \end{cases} \]
\end{theorem}

\begin{proof}
    Let $x_1$ and $x_2$ be the points in the hypotheses and set $R_1=|x_1|$ and $R_2=|x_2|$. Then $R_1>\bar r$ is the largest radius such that $B^{\star\star}(0,R_1)\subset\Omega^{\star\star}$, while $R_2>\bar r$ is the smallest radius such that $\Omega^{\star\star}\subset B^{\star\star}(0,R_2)$. Note that $R_1\le R_2$; see Figure~\ref{fig:conct_balls}.
\begin{figure}[h!]
    \centering
    \begin{tikzpicture}[scale=0.62,>=stealth]
        \draw[line width=0.8pt] (0,0) circle (0.85);
        \draw[line width=0.9pt] (0,0) circle (2.45);
        \draw[line width=0.9pt] (0,0) circle (4.00);

        \draw[BrickRed,line width=1.4pt]
          (-4.00,0)
          .. controls (-4.00,0.80) and (-3.70,1.50) .. (-2.75,2.35)
          .. controls (-1.55,3.35) and (0.60,3.70) .. (2.15,2.75)
          .. controls (3.20,2.10) and (3.55,0.75) .. (3.30,-0.55)
          .. controls (3.00,-1.75) and (1.45,-2.45) .. (0,-2.45)
          .. controls (-1.35,-2.45) and (-4.00,-0.80) .. (-4.00,0)
          -- cycle;

        \fill (0,0) circle (2.2pt);
        \fill (0,-2.45) circle (2.5pt);
        \fill (-4.00,0) circle (2.5pt);
        \node[font=\Large] at (-0.22,0.28) {$O$};
        \node[font=\Large] at (0,-3.00) {$X_1$};
        \node[font=\Large] at (-4.55,0.55) {$X_2$};

        \node[BrickRed,font=\Large] (Omega) at (3.85,3.70) {$\Omega^{**}$};
        \draw[BrickRed,->,line width=0.8pt]
          (Omega.south west) -- (2.65,2.35);

        \node[font=\large] (inner) at (1.70,1.35)
          {$B^{**}(0,\bar r)$};
        \draw[<-,line width=0.7pt] (inner.south west) -- (0.68,1.1);

        \node[font=\large] (Rone) at (4.15,-3.25)
          {$B^{**}(0,R_1)$};
        \draw[->,line width=0.7pt] (Rone.north west) -- (1.55,-1.90);

        \node[font=\large] (Rtwo) at (-4.15,-3.55)
          {$B^{**}(0,R_2)$};
        \draw[<-,line width=0.7pt] (Rtwo.north east) -- (-3.3,-3.1);
    \end{tikzpicture}
    \caption{The inner and outer concentric comparison balls and their contact points with $\partial\Omega^{\star\star}$.}
    \label{fig:conct_balls}
\end{figure}
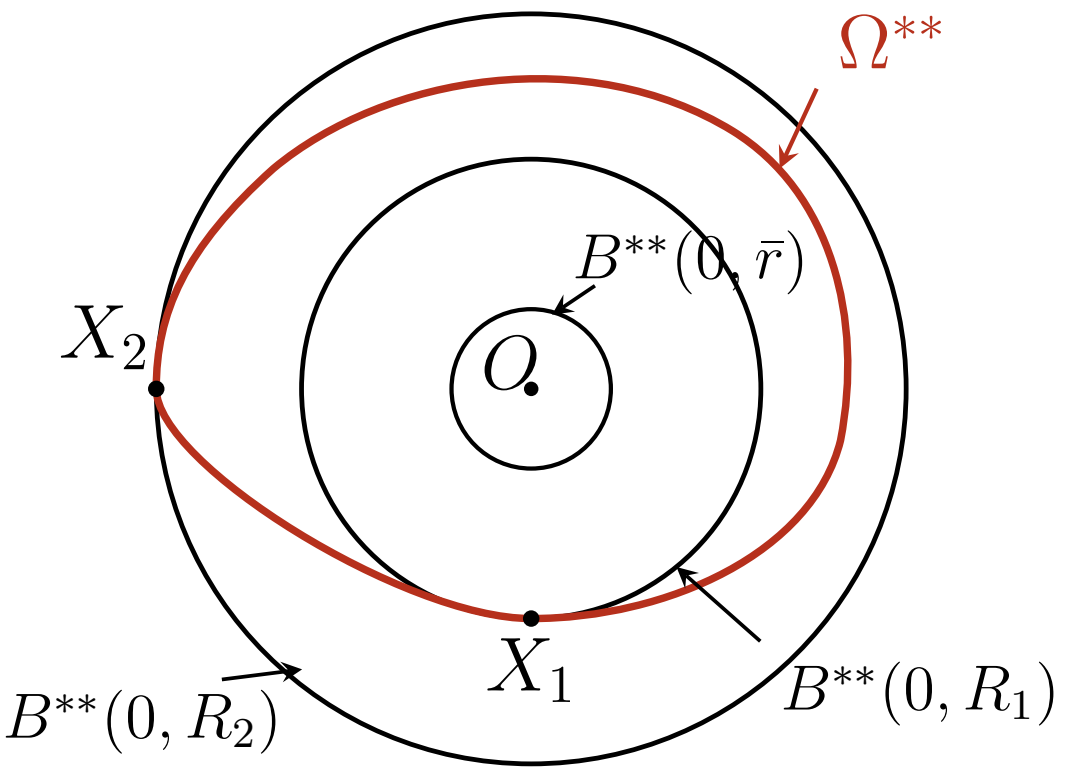
For $i=1,2$, let us define
    \begin{equation*}
    v_i(x)= \begin{cases}
     \displaystyle \frac{|x|^{\frac{p-n}{p-1}}- (R_i)^{\frac{p-n}{p-1}}}{\bar{r}^{\frac{p-n}{p-1}}- (R_i)^{\frac{p-n}{p-1}}} & p\neq n, \\
     \displaystyle  \frac{\ln(\frac{|x|}{R_i})}{\ln(\frac{\bar r}{R_i})}  & p= n. \\
    \end{cases}
\end{equation*}
    Then we have
    \[
    |\nabla  v_i |(x_i)=\begin{cases}
     \displaystyle  \frac{n-p}{p-1} \frac{(R_i)^{\frac{1-n}{p-1}}}{\bar{r}^{\frac{p-n}{p-1}}- (R_i)^{\frac{p-n}{p-1}}} & p< n, \\
     \displaystyle  \frac{p-n}{p-1} \frac{(R_i)^{\frac{1-n}{p-1}}}{(R_i)^{\frac{p-n}{p-1}} - \bar{r}^{\frac{p-n}{p-1}} } & p> n, \\
     \displaystyle   \frac{1}{R_i \ln(\frac{R_i}{\bar r})}  & p= n. \\
    \end{cases}
    \]
Here the comparisons are in the weak sense, since $G$ is not assumed continuous up to $\partial\Omega^{\star\star}$. Testing \eqref{eq:weakdirprob} with $G^-$ and $(G-1)^+$, which belong to $W^{1,p}_0(\Omega^{\star\star}\smallsetminus\overline{B^{\star\star}(0,\bar r)})$, gives $0\le G\le1$. Let $\tilde v_1=\min\{v_1,1\}$ in $B^{\star\star}(0,R_1)$. Then $0\le(\tilde v_1-G)^+\le\tilde v_1$ and $(\tilde v_1-G)^+=0$ in $B^{\star\star}(0,\bar r)$, so $(\tilde v_1-G)^+\in W^{1,p}_0(B^{\star\star}(0,R_1)\smallsetminus\overline{B^{\star\star}(0,\bar r)})$. Similarly $0\le(G-v_2)^+\le G$ and $(G-v_2)^+=0$ in $B^{\star\star}(0,\bar r)$, so $(G-v_2)^+\in W^{1,p}_0(\Omega^{\star\star}\smallsetminus\overline{B^{\star\star}(0,\bar r)})$. Testing the equations of $G$ and $v_i$ with these functions, the strict monotonicity of $\zeta\mapsto|\zeta|^{p-2}\zeta$ gives $G\ge v_1$ and $G\le v_2$ in the respective rings.

	    On one hand, since $G$ and $v_1$ are both $p$-harmonic in $B^{\star \star}(0, R_1) \smallsetminus \overline{B^{\star \star}(0,\bar r)}$ and $G \ge v_1$ on $\partial B^{\star \star}(0,\bar r) \cup \partial B^{\star \star}(0, R_1)$, we have $G\ge v_1$  in $B^{\star \star}(0, R_1) \smallsetminus \overline{B^{\star \star}(0,\bar r)}$, which implies
\[
|\nabla v_1(x_1)|=\limsup_{\substack{x \to x_1\\ x\in B^{\star\star}(0,R_1)}} \frac{v_1(x)}{|x-x_1|} \le
     \limsup_{x\to x_1} \frac{G(x)}{|x-x_1|}\le c.
    \]
    Hence we have
$$\begin{array}{cc}
   \frac{n-p}{p-1} \frac{(R_1)^{\frac{1-n}{p-1}}}{\bar{r}^{\frac{p-n}{p-1}}- (R_1)^{\frac{p-n}{p-1}}} \le c
,\quad \quad
&\textrm{for}\; p< n,\\
  \frac{p-n}{p-1} \frac{(R_1)^{\frac{1-n}{p-1}}}{(R_1)^{\frac{p-n}{p-1}} - \bar{r}^{\frac{p-n}{p-1}} } \le c
,
\quad \quad
&\textrm{for}\;  p> n,\\
      \frac{1}{R_1 \ln(\frac{R_1}{\bar r})} \le c,
&\textrm{for}\;   p= n.
\end{array}
$$
	    On the other hand, since $G$ and $v_2$ are both $p$-harmonic in $\Omega^{\star \star} \smallsetminus \overline{B^{\star \star}(0,\bar r)}$ and $G \le v_2$ on $\partial B^{\star \star}(0,\bar r)\cup  \partial \Omega^{\star \star}$, we have $G\le v_2$  in $\Omega^{\star \star} \smallsetminus \overline{B^{\star \star}(0,\bar r)}$, which gives
    \[
    |\nabla v_2(x_2)|
    \ge\limsup_{\substack{x\to x_2\\x\in\Omega^{\star\star}}}
    \frac{v_2(x)}{|x-x_2|}
    \ge\limsup_{\substack{x\to x_2\\x\in\Omega^{\star\star}}}
    \frac{G(x)}{|x-x_2|}\ge c.
    \]
 that is:
 $$
 \begin{array}{cc}
   \frac{n-p}{p-1} \frac{(R_2)^{\frac{1-n}{p-1}}}{\bar{r}^{\frac{p-n}{p-1}}- (R_2)^{\frac{p-n}{p-1}}} \ge c
,\quad \quad
&\textrm{for}\; p< n,\\
 \frac{p-n}{p-1} \frac{(R_2)^{\frac{1-n}{p-1}}}{(R_2)^{\frac{p-n}{p-1}} - \bar{r}^{\frac{p-n}{p-1}} } \ge c
,
\quad \quad
&\textrm{for}\;  p> n,\\
  \frac{1}{R_2\ln(\frac{R_2}{\bar r})} \ge c,
\quad \quad
&\textrm{for}\;   p= n.
\end{array}
$$
Then we can deduce
 $$
 \begin{array}{cc}
    \frac{(R_1)^{\frac{1-n}{p-1}}}{\bar{r}^{\frac{p-n}{p-1}}- (R_1)^{\frac{p-n}{p-1}}}   \le  \frac{(R_2)^{\frac{1-n}{p-1}}}{\bar{r}^{\frac{p-n}{p-1}}- (R_2)^{\frac{p-n}{p-1}}}
   ,\quad \quad
&\textrm{for}\; p< n,\\
     \frac{(R_1)^{\frac{1-n}{p-1}}}{(R_1)^{\frac{p-n}{p-1}} - \bar{r}^{\frac{p-n}{p-1}} }  \le   \frac{(R_2)^{\frac{1-n}{p-1}}}{(R_2)^{\frac{p-n}{p-1}} - \bar{r}^{\frac{p-n}{p-1}} }
,
\quad \quad
&\textrm{for}\;  p> n,\\
      \frac{1}{R_1 \ln(\frac{R_1}{\bar r})} \le  \frac{1}{R_2 \ln(\frac{R_2}{\bar r})}
,
\quad \quad
&\textrm{for}\;   p= n.
\end{array}
$$
Since the function $f:(\bar r, + \infty) \to \rr$ defined by \eqref{def:f}
	     is decreasing in $(\bar r, + \infty)$, it follows $R_2\le R_1$.

    Hence $R_2=R_1=R$ and $\Omega^{\star \star}=B^{\star \star}(0,R)$. Notice that $ R>0$ depends only on $c$ and $\bar r$ and it is the unique solution of the equations \eqref{eq:Rcr}, \eqref{eq:Rcr2}, \eqref{eq:Rcr3} (see the argument in the proof of Theorem \ref{th:euclideointerno}, replacing $2c$ with $c$).
    Moreover, by the uniqueness of the solution of the Dirichlet problem in $ B^{\star \star}(0,R) \smallsetminus B^{\star \star}(0,\bar r)$ we obtain that
    \[
    G(x)=v_1(x)=v_2(x)=\begin{cases} \vspace{0.3cm}
\displaystyle \frac{|x|^{\frac{p-n}{p-1}}- R^{\frac{p-n}{p-1}} }{\bar r^{\frac{p-n}{p-1}}- R^{\frac{p-n}{p-1}}} & p\ne n,
\\
\displaystyle \frac{\ln(\frac{|x|}{R})}{\ln(\frac{\bar r}{R})} & p=n.
 \end{cases}
    \]

\end{proof}
\begin{remark}
    The same argument of Remark \ref{rm:interiorBP} about the interior Bernoulli problem holds true also in this setting.
\end{remark}

We refer to the paragraph preceding Theorem \ref{th:subapprosimatointerno} for a reminder of Property $(H)$ (see Definition \ref{D:trans}), which is used in the next theorem. In addition, we shall use \eqref{e:circimscribed and inscribed balls} to define the points  of contact $\operatorname{Min}({\partial\Omega})$ and $
\operatorname{Max}({\partial\Omega})
$ of the inner and outer gauge balls centered at the origin that contain $\partial \Omega$. Notice that when  $\Om$ has partial symmetry, since the gauge has partial symmetry (see Definition \ref{df:partialsy}), the sets of contact of  $\partial\Om$ and the outer and inner balls also have partial symmetry, i.e., they are invariant under the sub-group  $\mathscr S(\frak{z}^{\perp})$  of orthogonal maps of $\frak{z}^{\perp}$ which commute with all of the maps $J_{\sigma}$ in \eqref{e:Jz def}, $\sigma$ a unit length vector in $\frak{z}$.

\begin{theorem}
\label{th:subapplicationtosharp}
Suppose that $\bG$ satisfies Property $(H)$ and let $p \in \nn$ be an even number. Let $\Omega \subset \bG$ be a connected, bounded open set  with partial symmetry such that $0 \in \Omega$, and there exists $r>0$ such that $B_r(0) \Subset \Omega$. Let $u$ be a weak solution of
\begin{equation}
\label{eq:Beprb2}
 \begin{cases}
        \Delta_{\bG,p} u=0 &  \Omega \smallsetminus  \overline{B_r(0)},\\
        u=0 & \partial \Omega,\\
        u=1 &\overline{ B_r(0)}.\\

    \end{cases}
\end{equation}

Assume that for some $c>0$ there exist $(z_1,\sigma_1) \in \operatorname{Min}({\partial\Omega}) $ and $(z_2,\sigma_2) \in \operatorname{Max}({\partial\Omega}) $, cf. \eqref{e:circimscribed and inscribed balls}, with respect to the norm $|\cdot|_{\bG}$, such that $u(z,\sigma)\to0$ as $(z,\sigma)\to(z_1,\sigma_1)$, $(z,\sigma)\in\Omega$, and the following $\limsup$ is taken only over points with $z\ne0$:
\begin{equation}
\label{eq:limsuples}
\limsup_{\substack{(z,\sigma)\to(z_1,\sigma_1)\\(z,\sigma)\in\Omega,\ z\ne0}}
\frac{|\nabla_{\bG}u(z,\sigma)|}{|z|}\le c,
\end{equation}
and assume that, if $z_2 \neq 0$, there exists $t^* \in (0, 1)$ such that $(t z_2,\sigma_2) \in \Omega$ for all $t \in (t^*,1)$ and
\begin{equation}
\label{subcond:exterior_sharp}
\limsup_{t  \to 1^-} \frac{u (t z_2,\sigma_2)}{(1-t)|z_2|^2} \ge {c};
\end{equation}
and, if $z_2=0$, there exists $t^* \in (0, 1)$ such that $(0,t\sigma_2) \in \Omega$ for all $t \in (t^*,1)$ and
\begin{equation}
\label{subcond:exterior2}
\limsup_{t  \to 1^-} \frac{u (0,t\sigma_2)}{ (1-t)|\sigma_2|} \ge 2c.
\end{equation}

Then, there exists $R>r$ such that $\Omega=B_R(0)$ and
\begin{equation}
\label{eq:usol2}
u(z,\sigma)=\begin{cases}
\vspace{0.3 cm}
\displaystyle \frac{|(z, \sigma)|_{\bG}^{\frac{p- Q}{p-1}} - R^{\frac{p- Q}{p-1}}}{r^{\frac{p- Q}{p-1}}- R^{\frac{p- Q}{p-1}}} & p \neq Q,\\
\displaystyle \frac{\ln(\frac{|(z,\sigma)|_{\bG}}{R})}{\ln(\frac{r}{R}) } & p=Q.
\end{cases}
\end{equation}
\end{theorem}

\begin{proof}
As in the proof of Theorem \ref{th:subapprosimatointerno}, $u$ is partially symmetric and, with $n$, $\bar r$, $\Omega^{\star\star}$, $G$, $\Pi$, $\Psi$ as there:
\begin{itemize}
\item $\Omega^{\star\star}$ is a connected bounded open set with $\overline{B^{\star\star}(0,\bar r)}\subset\Omega^{\star\star}$;
\item $G\in W^{1,p}_0(\Omega^{\star\star})$ satisfies \eqref{eq:weakdirprob};
\item $G=1$ a.e. in $\overline{B^{\star\star}(0,\bar r)}$;
\item $G\in C^{1,\alpha}_{loc}(\Omega^{\star\star}\smallsetminus\overline{B^{\star\star}(0,\bar r)})$, and \eqref{eq:gradcorr} holds.
\end{itemize}
 Since $\Pi(\partial\Omega)\subset\partial\Omega^{\star\star}$, $\Psi(\partial\Omega^{\star\star})\subset\partial\Omega$, $|\Pi(y)|=\frac{|y|_{\bG}^2}{4}$ and $|\Psi(x)|_{\bG}^2=4|x|$, the points $x_1=\Pi(z_1,\sigma_1)$ and $x_2=\Pi(z_2,\sigma_2)$ belong to $\operatorname{Min}(\partial\Omega^{\star\star})$ and $\operatorname{Max}(\partial\Omega^{\star\star})$ respectively. We verify the hypotheses of Theorem \ref{th:eucsharp} with $2c$ in place of $c$.

Let $x=(w,\sigma)\in\Omega^{\star\star}$ tend to $x_1$. Let $\zeta=z_1/|z_1|$ if $z_1\ne0$, $\zeta=e$ otherwise, and $y=(2\sqrt{|w|}\,\zeta,\sigma)$. Then $y\in\Omega$, $y\to(z_1,\sigma_1)$, $G(x)=u(y)$, and, by \eqref{eq:gradcorr}, $|\nabla G(x)|=2|\nabla_{\bG}u(y)|/(2\sqrt{|w|})$ if $w\ne0$ (recall that $G(w,\sigma)=g(|w|,\sigma)$, so that $G$ and $|\nabla G|$ are invariant under rotations in $w$). By \eqref{eq:limsuples} and the continuity of $\nabla G$,
\[
\limsup_{\substack{x\to x_1\\x\in\Omega^{\star\star}}}|\nabla G(x)|\le2c,
\]
and, by the hypothesis $u\to0$ at $(z_1,\sigma_1)$, $G(x)\to0$ as $x\to x_1$, $x\in\Omega^{\star\star}$. Fix $\eps>0$ and $0<\rho<\operatorname{dist}(x_1,\overline{B^{\star\star}(0,\bar r)})$ with $|\nabla G|\le2c+\eps$ in $\Omega^{\star\star}\cap B_\rho(x_1)$. The extension $\widetilde G$ of $G$ by zero belongs to $W^{1,p}(\rr^n)$, with $\nabla\widetilde G=\mathbf 1_{\Omega^{\star\star}}\nabla G$ a.e. Hence $|\nabla\widetilde G|\le2c+\eps$ a.e. in $B_\rho(x_1)$, and $\widetilde G$ has a $(2c+\eps)$-Lipschitz representative in $B_\rho(x_1)$. This representative coincides with $G$ in $\Omega^{\star\star}\cap B_\rho(x_1)$, both being continuous there, and therefore vanishes at $x_1$. Thus $G(x)\le(2c+\eps)|x-x_1|$ in $\Omega^{\star\star}\cap B_\rho(x_1)$, and, $\eps$ being arbitrary,
\[
\limsup_{\substack{x\to x_1\\x\in\Omega^{\star\star}}}\frac{G(x)}{|x-x_1|}\le2c .
\]
No segment from $x$ to $x_1$ is required to lie in $\Omega^{\star\star}$.

If $z_2\ne0$, then $x_2=(w_2,\sigma_2)$ with $|w_2|=\frac{|z_2|^2}{4}$, and for $t\in(t^*,1)$ the points $(t^2w_2,\sigma_2)=\Pi(tz_2,\sigma_2)$ belong to $\Omega^{\star\star}$, with $G(t^2w_2,\sigma_2)=u(tz_2,\sigma_2)$ and $|(t^2w_2,\sigma_2)-x_2|=(1-t^2)\frac{|z_2|^2}{4}$. By \eqref{subcond:exterior_sharp},
\[
\limsup_{\substack{x\to x_2\\x\in\Omega^{\star\star}}}\frac{G(x)}{|x-x_2|}\ge\limsup_{t\to1^-}\frac{4\,u(tz_2,\sigma_2)}{(1+t)(1-t)|z_2|^2}=2\limsup_{t\to1^-}\frac{u(tz_2,\sigma_2)}{(1-t)|z_2|^2}\ge2c .
\]
If $z_2=0$, then $x_2=(0,\sigma_2)$, the points $(0,t\sigma_2)=\Pi(0,t\sigma_2)$, $t\in(t^*,1)$, belong to $\Omega^{\star\star}$, and by \eqref{subcond:exterior2}
\[
\limsup_{\substack{x\to x_2\\x\in\Omega^{\star\star}}}\frac{G(x)}{|x-x_2|}\ge\limsup_{t\to1^-}\frac{u(0,t\sigma_2)}{(1-t)|\sigma_2|}\ge2c .
\]
By Theorem \ref{th:eucsharp}, applied with $2c$ in place of $c$, there is $\bar R>\bar r$ such that $\Omega^{\star\star}=B^{\star\star}(0,\bar R)$ and $G$ is the radial function given there. As at the end of the proof of Theorem \ref{th:subapprosimatointerno}, $\Omega=B_R(0)$ with $R=2\bar R^{\frac12}$, and $u=G\circ\Pi$ is given by \eqref{eq:usol2}.

\end{proof}

\end{document}